\documentclass[10pt,final]{amsart}

\usepackage{graphicx}
\usepackage{amsmath}
\usepackage{amssymb}
\usepackage[latin1]{inputenc}
\usepackage[all]{xy}
\usepackage{semantic}
\usepackage{rotating}
\usepackage{hyperref}
\usepackage{enumitem}
\usepackage{relsize}
\usepackage{booktabs}
\usepackage{color}
\usepackage{tikz}
\usepackage{tikz-cd}

\theoremstyle{plain} \newtheorem{theorem}{Theorem}[section]
\theoremstyle{plain} \newtheorem{proposition}[theorem]{Proposition}
\theoremstyle{plain} \newtheorem{corollary}[theorem]{Corollary}
\theoremstyle{plain} \newtheorem{lemma}[theorem]{Lemma}
\theoremstyle{definition} \newtheorem{definition}[theorem]{Definition}
\theoremstyle{definition} \newtheorem{remark}[theorem]{Remark}
\theoremstyle{definition} \newtheorem{example}[theorem]{Example}

\newcommand{\cS}{\mathcal{S}}

\renewcommand{\phi}{\varphi}

\newcommand{\ua}{\ensuremath{\underline a}}
\newcommand{\ub}{\ensuremath{\underline b}}
\newcommand{\uc}{\ensuremath{\underline c}}

\newcommand{\ux}{\ensuremath{\underline x}}

\newcommand{\uy}{\ensuremath{\underline y}}
\newcommand{\uz}{\ensuremath{\underline z}}
\newcommand{\up}{\ensuremath{\underline p}}

\newcommand{\uw}{\ensuremath{\underline w}}

\newcommand{\cA}{\ensuremath \mathcal A}
\newcommand{\cB}{\ensuremath \mathcal B}
\newcommand{\cM}{\ensuremath \mathcal M}
\newcommand{\cN}{\ensuremath \mathcal N}

\newcommand{\cC}{\ensuremath \mathcal C}

\newcommand{\upsi}{\ensuremath{\underline \psi}}
\newcommand{\uphi}{\ensuremath{\underline{\varphi}}}
\renewcommand{\cong}{\wedge}

\newcommand{\dis}{\vee}
\renewcommand{\ne}{\neg}
\newcommand{\impl}{\rightarrow}
\newcommand{\biimpl}{\leftrightarrow}
\renewcommand{\tt}{\top}
\newcommand{\ff}{\bot}

\newcommand{\univ}{[\forall]}
\newcommand{\simpa}{\leadsto}
\newenvironment{ourlist}{\begin{list}{}
    {\setlength{\topsep}{1mm}\setlength{\itemsep}{0mm}\setlength{\parsep}{0mm}}
    }{\end{list}}

\newcommand{\amp}{{ \ \& \ }}
\newcommand{\relr}{\mathrel{R}}
\def\Pcal{\mathcal{P}}
\def\Str{{\sf{StR}}}
\def\Con{{\sf{Con}}}
\def\S{\mathcal{S}}
\def\fowedge{\amp}
\def\fovee{ \ \mathtt{or} \ }
\def\foneg{{ \sim }} 

\usepackage{fullpage} 

\title{{Admissibility of $\Pi_2$-Inference Rules: interpolation, model completion, and contact algebras}}
\date{}
\author{Nick Bezhanishvili}
\author{Luca Carai}
\author{Silvio Ghilardi}
\author{Lucia Landi}

\begin{document}

\begin{abstract} 
We devise three strategies for recognizing admissibility of non-standard inference rules via interpolation, uniform interpolation, and model completions. 
We apply our machinery to the case of symmetric implication calculus $\mathsf{S^2IC}$, where we also supply a finite axiomatization of the model completion
of its algebraic counterpart, via the equivalent theory of contact algebras. Using this result we obtain a finite basis for admissible $\Pi_2$-rules. 
\end{abstract}


\maketitle 
 
\section{Introduction}

The use of non-standard rules has a long tradition in modal logic starting from the pioneering work of  Gabbay \cite{Gab81}, who
introduced a non-standard rule for irreflexivity.  Non-standard rules have been employed in temporal logic   in
the context of branching time logic \cite{Bur80}  and for axiomatization problems \cite{GH90} concerning 
the logic of the real line in the language with the Since and Until modalities. General completeness
results for modal languages that are sufficiently expressive to define the so-called difference modality have been obtained in  \cite{Ven93}. 
For the use of the non-standard density rule in many-valued logics we refer to  \cite{MM07} and \cite{TT84}.

Recently, there has been a renewed interest in non-standard rules in the context of the region-based theories of space \cite{Vak07}. One of the key algebraic structures 
in these theories is that of \emph{contact algebras}. These algebras form a discriminator variety, see, e.g., \cite{BBSV}. Compingent algebras are contact algebras satisfying 
two $\forall\exists$-sentences (aka $\Pi_2$-sentences)  \cite{BBSV,DeV62}. De Vries  \cite{DeV62} established a duality between complete compingent algebras and compact Hausdorff spaces. 
This duality led to new logical calculi for 
compact Hausdorff spaces in \cite{BTV07} for a two-sorted modal language and in \cite{BBSV} for a uni-modal language with a strict implication. Key to these approaches is a development  of  logical 
calculi corresponding to contact algebras. In \cite{BBSV} such a calculus is called the \emph{strict symmetric implication calculus} and is denoted by  $\mathsf{S^2IC}$. The extra $\Pi_2$-axioms of compingent algebras then correspond to non-standard $\Pi_2$-rules, which turn out to be admissible in  $\mathsf{S^2IC}$. This generates a natural question of investigating admissibility of $\Pi_2$-rules in $\mathsf{S^2IC}$ studied in~\cite{BBSV} and in general in logical calculi corresponding to  
varieties of modal algebras. In fact, rather little is known about the problem of recognizing \emph{admissibility} for such  non-standard rules, although this problem has already been raised in \cite{Ven93}.
This is the question that we address in this paper. 

We undertake a systematic study of admissibility of $\Pi_2$-rules. As far as we are aware, this is a first
attempt to study admissibility in the context of non-standard inference rules. 
In fact, we show that \emph{there are  tools already available in the literature on modal logic} that can be fruitfully employed for this aim: these tools include   algorithms for deciding conservativity, as well as algorithms for computing local and global interpolants. 
We devise three different strategies for recognizing admissibility of $\Pi_2$-rules over some system $\mathcal{S}$. The definition of $\Pi_2$-rules that we consider is taken from \cite{BBSV} and is close
to that of Balbiani et al.~\cite{BTV07}.

The first strategy applies to a logic $\mathcal{S}$ with the interpolation property.
 We show that $\Pi_2$-rules are effectively recognizable in $\mathcal{S}$ in case $\mathcal{S}$ has 
the interpolation property and conservativity is decidable in $\mathcal{S}$.  
The second strategy applies  to logics admitting local and global uniform interpolants, respectively.
Global interpolants are strictly related to model completions and to axiomatizations of existentially closed structures~\cite{GZ}, thus establishing a direct connection between $\Pi_2$-rules and model-theoretic machinery. 
Directly exploiting this connection leads to our third strategy.
We apply the third strategy to our main case study 
to show admissibility of various $\Pi_2$-rules in $\mathsf{S^2IC}$, thus recovering admissibility results from~\cite{BBSV} as special cases (we also show that the admissibility problem  for $\mathsf{S^2IC}$ is co-\textsc{NExpTime}-complete). 
The model completion we use to this aim is that of the theory of contact algebras.

Finally, we prove  the technically most challenging result of the paper: that the model completion of contact algebras is  finitely axiomatizable. As a consequence of this result  we obtain a finite basis for admissible $\Pi_2$-rules in $\mathsf{S^2IC}$.

\section{Preliminaries}\label{sec:preliminaries}

A \emph{modal signature} $\Sigma$ is a finite signature comprising Boolean operators $\cong, \dis, \impl,\biimpl, \ne$ as well as additional operators of any arity called the \emph{modal} operators. 
Out of $\Sigma$-symbols and out of a countable set of variables  $x, y, z, \dots, p,q,r,\dots$ one can build the set of propositional \emph{$\Sigma$-formulas}. 
$\Sigma$-formulas might be indicated both with the greek letters $\phi, \psi, \dots$ and the latin capital letters $F, G, \dots$. Notations such as  $F(\ux)$ mean that the $\Sigma$-formula $F$ contains at most the variables from the tuple $\ux$; the notation $F(\underline{\varphi}/\ux)$ means the simultaneous componentwise substitution of the tuple of formulas  $\underline{\varphi}$ for the tuple of variables $\ux$. 

A \emph{modal system $\cS$} (over the modal signature $\Sigma$) is a set of $\Sigma$-formulas comprising tautologies, the distribution axioms~\footnote{ 
Extension to  non normal operators needs to be investigated.
} 
$$
\Box[\phi, \dots, \psi\to \psi',\dots, \phi]\to
(\Box[\phi, \dots, \psi,\dots, \phi] \to
\Box[\phi, \dots, \psi',\dots, \phi])
$$
and
closed under the rules of modus ponens (MP) 
(from $\phi$ and $\phi\impl\psi$ infer $\psi$), uniform substitution (US) (from $F(\ux)$ infer $F(\underline{\psi}/\ux)$),  
and necessitation  (N) 
(from $\psi$  infer
$\Box[\phi, \dots, \psi,\dots, \phi]$).

We write $\vdash_{\cS} \phi$ or $\cS\vdash \phi$ to mean that $\phi\in \cS$. If $\vdash_{\cS} \phi\to \psi$ holds, we say that $\psi$ is a \emph{local} consequence of $\phi$ (modulo $\cS$). We shall also need the \emph{global} consequence relation $\phi\vdash_{\cS} \psi$: this relation holds when $\psi$ belongs to the smallest set of formulas containing $\cS$ and  $\phi$ that is closed under modus ponens and 
necessitation
(notice that closure under uniform substitution is not required).

We say that $\cS$ is decidable iff the relation $\vdash_\cS \phi$ is decidable. We also say that $\cS$ is \emph{locally tabular} iff for every finite tuple of propositional variables $\ux$ there are finitely many formulas $\psi_1(\ux), \dots, \psi_n(\ux)$ such that for every further formula $\phi(\ux)$ there is some $i=1, \dots, n$ such that $\vdash_\cS \phi\biimpl \psi_i$.

We say that $\cS$ has the (local) \emph{interpolation property} iff for every pair of $\Sigma$-formulas $\phi(\ux, \uy), \psi(\uy, \uz)$ such that 
$\vdash_{\cS} \phi\to \psi$ there is a formula $\theta(\uy)$ such that $\vdash_{\cS} \phi\to \theta$ and $\vdash_{\cS} \theta\to \psi$. Similarly, we say that $\cS$ has the  \emph{global interpolation property} iff for every pair of $\Sigma$-formulas $\phi(\ux, \uy), \psi(\uy, \uz)$ such that 
$\phi\vdash_{\cS}  \psi$ there is a formula $\theta(\uy)$ such that $\phi\vdash_{\cS} \theta$ and $\theta \vdash_{\cS}  \psi$. Using Lemma~\ref{lem:ded} below it is easily seen that the local interpolation property implies the global interpolation property (the converse however does not hold, even over $\mathsf{S4}$, see~\cite{Mak,Mak1}).

Let us call a \emph{definable modality} or simply a \emph{modality} any formula $M(x)$ (where only the variable $x$ occurs) such that $x\vdash_{\cS} M(x)$ and $\vdash_{\cS} M(x_1\wedge x_2/x) \leftrightarrow  M(x_1/x)\wedge M(x_2/x)$ (notice that $\vdash_{\cS} M(\top/x)$ follows as a consequence). 
 
\begin{lemma} [Deduction]\label{lem:ded} If $\phi\vdash_{\cS} \psi$ holds, there is a modality $M(x)$ such that 
 $\vdash_{\cS} M(\phi/x)\to \psi$.
\end{lemma}
\begin{proof}
The required modality \emph{depends on the derivation} and is built up inductively as follows.
For length 1 derivations consisting of axioms or of $\phi$, we take $M(x)$ to be $x$.
If $\psi$ is obtained from $\psi'\to \psi$ and $\psi'$ via modus ponens, then by induction we have modalities $M_1(x), M_2(x)$ such that 
$\vdash_{\cS} M_1(\phi/x)\to (\psi'\to \psi)$ and $\vdash_{\cS} M_2(\phi/x)\to \psi'$.
Then we get $\vdash_{\cS} M_1(\phi/x)\wedge M_2(\phi/x)\to \psi$ and $M_1(x)\wedge M_2(x)$ is our desired modality.
 If $\psi$ is obtained by the necessitation rule, it is of the kind 
 $\Box[\theta, \dots, \psi', \dots,\theta]$
 and we have by induction a modality $M(x)$ such that 
 $\vdash_{\cS} M(\phi/x)\to \psi'$; we then get by necessitation and normality
 $$
 \vdash_{\cS}\Box[\theta, \dots, M(\phi/x), \dots,\theta]\to
 \Box[\theta, \dots, \psi', \dots,\theta].
 $$
 Since $\vdash_{\cS} \bot \to \theta$, by iterated applications of necessitation and normality, we  obtain
 $$
 \vdash_{\cS}
 \Box[\bot, \dots, M(\phi/x), \dots,\bot]
 \to \Box[\theta, \dots, M(\phi/x), \dots,\theta]
 $$
 and also 
 $$
 \vdash_{\cS}\Box[\bot, \dots, M(\phi/x), \dots,\bot]\to
 \Box[\theta, \dots, \psi', \dots,\theta]
 $$
 by transitivity of implication. 
Thus, $\vdash_{\cS}\Box[\bot, \dots, M(\phi/x), \dots,\bot]\to \psi$.
It is straightforward to see that $\Box[\bot, \dots, M(x), \dots,\bot]$ is a modality.
\end{proof}

\begin{lemma}[Replacement]\label{lem:repl}
 For every $n$-tuple of variables $\ux:=x_1, \dots, x_n$, for every  formula $\phi(\ux, \uy)$, 
 and for every pair of $n$-tuples of formulas $\underline{\psi}:=\psi_1,\dots, \psi_n$, $\underline{\psi}':=\psi'_1, \dots \psi'_n$ we have that 
 $$
 \bigwedge_{i=1}^n \psi_i\leftrightarrow \psi_i' \vdash_{\cS}
 \phi(\underline{\psi}/\ux,\uy)\leftrightarrow
 \phi(\underline{\psi}'/\ux, \uy).
 $$
 As a consequence, by Lemma~\ref{lem:ded}, there is a modality $M(x)$ (depending on $\phi,\underline{\psi},\underline{\psi}'$)  such that 
 $$
 \vdash_{\cS} \bigwedge_{i=1}^n M(\psi_i\leftrightarrow \psi_i'/x) \wedge \phi(\underline{\psi}/\ux,\uy)\to
 \phi(\underline{\psi}'/\ux, \uy).
 $$
\end{lemma}
\begin{proof}
We prove the statement by induction on $\phi$. If $\phi$ is a propositional variable or its main connective is a Boolean connective, the statement is trivial. If the main connective of $\phi$ is a Box operator, it is sufficient to see that the following replacement rule 
$$
\rm{from~} \psi\leftrightarrow \psi'~{\rm infer}~ 
\Box[\theta, \dots, \psi, \dots,\theta] \leftrightarrow
\Box[\theta, \dots, \psi', \dots,\theta]
$$
is derivable. In fact, if $\vdash_{\cS} \psi\leftrightarrow \psi'$, then 
$\vdash_{\cS} \psi\rightarrow \psi'$ and $\vdash_{\cS}
\Box[\theta, \dots, \psi, \dots,\theta] \rightarrow
\Box[\theta, \dots, \psi', \dots,\theta]$ follow. Analogously,  
$\Box[\theta, \dots, \psi', \dots,\theta] \rightarrow
\Box[\theta, \dots, \psi, \dots,\theta]$ and finally
$\Box[\theta, \dots, \psi, \dots,\theta] \leftrightarrow
\Box[\theta, \dots, \psi', \dots,\theta]$.
\end{proof}

We say that $\cS$ has \emph{universal modality} iff $\Sigma$ contains a unary operator $\univ$ and $\cS$ includes the following formulas:
\begin{eqnarray*}
&\univ \phi \impl \phi,
&\univ \phi \impl \univ \univ \phi, \\
&\phi \impl \univ \ne\univ\ne \phi,
&\univ(\phi\impl \psi)\impl (\univ \phi \impl \univ \psi),\\
&\bigwedge_i \univ (\phi_i \biimpl \psi_i) \impl (\Box[\phi_1, \dots, \phi_n]\biimpl \Box[\psi_1, \dots, \psi_n]) ~ &({\rm for~all}~\Box\in \Sigma).
\end{eqnarray*}
For systems with universal modalities, Lemmas~\ref{lem:ded} and~\ref{lem:repl} can be simplified as follows:
\begin{lemma} [Deduction-Replacement]\label{lem:dedrepl}
 Let $\cS$ have a universal modality; then
 \begin{description}
  \item[{\rm (i)}] $\phi\vdash_{\cS} \psi$ holds iff 
 $\vdash_{\cS} \univ \phi\to \psi$;
 \item[{\rm (ii)}] the following sentences are provable
 $$
  \bigwedge_{i=1}^n \univ(\psi_i\leftrightarrow \psi_i') \wedge \phi(\underline{\psi}/\ux,\uy)\to
 \phi(\underline{\psi}'/\ux, \uy)
 $$
 for every pair of $n$-tuples of formulas $\underline{\psi}:=\psi_1,\dots, \psi_n$, $\underline{\psi}':=\psi'_1, \dots, \psi'_n$.
 \end{description}

\end{lemma}

We finally introduce $\Pi_2$-rules, which are the main objects of study of this paper.

\begin{definition} \label{defirules}
A \emph{$\Pi_2$-rule} is a rule of the form
\[
(\rho) \quad \inference{F(\underline{\varphi}/\ux,\up) \impl \chi}
   {G(\underline{\varphi}/\ux) \impl \chi}
\]
where $F(\ux, \up),G(\ux)$ are formulas.
We say that $\theta$ is obtained from $\psi$ by an application
of the rule $\rho$ if $\psi = F(\underline{\phi}/\ux,\up) \to \chi$ and
$\theta = G(\underline{\phi}/\ux) \to \chi$, where $\underline{\varphi}$ is a tuple of formulas, $\chi$ is a formula,
and $\up$ is a tuple of 
propositional letters not occurring in $\underline{\varphi}, \chi$.
\end{definition}

The definition of $\Pi_2$-rules follows the one of \cite{BBSV} and is close
to that of Balbiani et al.~\cite{BTV07}.
We now consider the effect of the addition of $\Pi_2$-rules to a system $\cS$.

\begin{definition}[Proofs with $\Pi_2$-rules] \label{proofrules}
Let $\Theta$ be a set of $\Pi_2$-rules. For 
a formula $\varphi$, we say that $\varphi$ is \emph{derivable}
in $\cS$ using the $\Pi_2$-rules in $\Theta$, and write 
$\vdash_{\cS+\Theta}\varphi$, provided
there is a proof $\psi_1,\dots,\psi_n$ such that $\psi_n=\varphi$ and each $\psi_i$ 
is an instance of an axiom of $\cS$, or is
obtained either by (MP), (N) or by an application of a rule $\rho \in \Theta$ from some previous $\psi_j$'s.
\end{definition}

We are interested in characterizing those $\Pi_2$-rules that can be freely used in a system without affecting its deductive power.

\begin{definition}
A rule $\rho$ is \emph{admissible} in the system $\mathcal S$ if for each formula $\varphi$, from $\vdash_{\mathcal{S}+\rho} \varphi$
it follows that $\vdash_{\mathcal S} \varphi$. The \emph{admissibility problem for $\Pi_2$-rules} is the following: given a $\Pi_2$-rule $\rho$ and a system $\cS$, decide whether it is admissible or not in $\cS$.
\end{definition}

In the rest of the paper we will study admissibility of $\Pi_2$-rules. 
\section{Conservative Extensions}\label{sec:conservative}

In this section we describe how to determine whether a $\Pi_2$-rule is admissible via conservative extensions.
Conservative extensions in modal logics were  investigated in~\cite{GLWZ06} 
and in description logics in~\cite{consbib0,consbib1,consbib2};
we recall here the related definition:

\begin{definition}
 Let $\phi_1(\ux), \phi_2(\ux, \uy)$ be $\Sigma$-formulas; we say that $\phi_1(\ux)\wedge \phi_2(\ux,\uy)$ is a conservative extension of $\phi_1(\ux)$ in $\cS$ iff for every further $\Sigma$-formula $\psi(\ux)$, we have that
 $$
 \vdash_{\cS}\phi_1\wedge \phi_2 \to \psi \quad \Rightarrow \quad
 \vdash_{\cS}\phi_1\to \psi.
 $$
 The \emph{conservativity problem for $\cS$} is the following: given $\phi_1(\ux)$ and $ \phi_2(\ux, \uy)$ decide whether that $\phi_1(\ux)\wedge \phi_2(\ux,\uy)$ is a conservative extension of $\phi_1(\ux)$ in $\cS$.
\end{definition}

\begin{theorem}\label{thm:cons}
 Assume that $\cS$ has the interpolation property. 
 Then a $\Pi_2$-rule $\rho$ of the form
 \[
\inference{F(\underline{\varphi}/\ux,\up) \impl \chi}
   {G(\underline{\varphi}/\ux) \impl \chi}
\]
 is admissible in $\cS$ iff 
 $G(\ux)\wedge F(\ux, \up)$ is a conservative extension  of $G(\ux)$ in $\cS$.
\end{theorem}

\begin{proof}
 For the left-to-right side, assume that $\rho$ is admissible and that  $\vdash_{\cS} F(\ux, \up) \wedge G(\ux) \to H(\ux)$. Then we have $\vdash_{\cS} F(\ux, \up) \to ( G(\ux) \to H(\ux))$,  and by admissibility $\vdash_{\cS} G(\ux) \to ( G(\ux) \to H(\ux))$ which is the same as  
 $\vdash_{\cS}  G(\ux) \to H(\ux)$. This shows that $G(\ux)\wedge F(\ux, \up)$ is a conservative extension  of $G(\ux)$ in $\cS$.
 
 For the converse, assume that $G(\ux)\wedge F(\ux, \up)$ is a conservative extension  of $G(\ux)$ in $\cS$ and that 
 \begin{equation*}
  \vdash_{\cS} F(\underline{\phi}(\uy)/\ux, \up)  \to \chi(\uy)~.
 \end{equation*}
 Let $\ux=x_1, \dots, x_n$; 
 since we have 
 $$
 \bigwedge_{i=1}^n (x_i\leftrightarrow \phi_i(\uy)) \vdash_{\cS}  F(\ux, \up)  \leftrightarrow 
  F(\underline{\phi}(\uy)/\ux, \up),
 $$
 by applying the Replacement Lemma~\ref{lem:repl}, we obtain a modality $M(x)$
such that
 \begin{equation*}
  \vdash_{\cS} \bigwedge_{i=1}^n M(x_i\leftrightarrow \phi_i(\uy)) \wedge F(\ux, \up)  \to 
  F(\underline{\phi}(\uy)/\ux, \up).
 \end{equation*}
Thus, by transitivity of implication, we have
 \begin{equation*}
  \vdash_{\cS} \bigwedge_{i=1}^n M(x_i\leftrightarrow \phi_i(\uy)) \wedge F(\ux, \up)  \to \chi(\uy);
 \end{equation*}
that is equivalent to
 \begin{equation*}
  \vdash_{\cS}  F(\ux, \up)  \to \left( \bigwedge_{i=1}^n M(x_i\leftrightarrow \phi_i(\uy)) \to \chi(\uy) \right).
 \end{equation*}
 By the interpolation property, there is $\theta(\ux)$ such that
 \begin{equation}\label{eq:from_interp}
  \vdash_{\cS}  F(\ux, \up)  \to \theta(\ux)
  ~{\rm and}~ \vdash_{\cS} \theta(\ux)\to \left( \bigwedge_{i=1}^n M(x_i\leftrightarrow \phi_i(\uy)) \to \chi(\uy) \right).
 \end{equation}
 The former entailment implies $\vdash_{\cS} G(\ux)\wedge  F(\ux, \up)  \to \theta(\ux)$ and so, by conservativity, we get $\vdash_{\cS} G(\ux) \to \theta(\ux)$.
 From the second entailment of~\eqref{eq:from_interp}, by transitivity, we then obtain
 \begin{equation*}
  \vdash_{\cS} G(\ux)\to \left( \bigwedge_{i=1}^n M(x_i\leftrightarrow \phi_i(\uy)) \to \chi(\uy) \right).
 \end{equation*}
Applying the replacements $\phi_i(\uy)/x_i$, we finally get 
 \begin{equation*}
  \vdash_{\cS} G(\underline{\phi}(\uy)/\ux)\to  \chi(\uy),
 \end{equation*}
 showing that $\rho$ is admissible.
\end{proof}

\begin{remark}
We underline that, 
without interpolation, the right-to-left implication of Theorem~\ref{thm:cons} may fail. In fact, let $\cS$ be a locally tabular logic without interpolation (see~\cite{Mak, Mak1} for examples) and suppose that 
interpolation fails for the entailment $\vdash_{\cS} F(\ux,\up)\to H(\ux,\uw)$. Let $G(\ux)$ be the conjunction of the $\ux$-formulas implied 
by $F$ (up to logical equivalence, there are only finitely many of them). Then $G \wedge F$ is obviously conservative over $G$. However, the relative $\Pi_2$-rule 
\[
\inference{F(\underline{\varphi}/\ux,\up) \impl \chi}
   {G(\underline{\varphi}/\ux) \impl \chi}
\]
is not admissible. Indeed, if it were, from $\vdash_{\cS} F(\ux,\up)\to H(\ux,\uw)$
 we would obtain $\vdash_{\cS} G(\ux)\to H(\ux,\uw)$ implying that  $G(\ux)$ is  an interpolant.
\end{remark}

Thus, we have obtained the following:

\begin{corollary}
If $\cS$ has the interpolation property and conservativity is decidable in $\cS$, then $\Pi_2$-rules are effectively recognizable in $\cS$. 
\end{corollary}

\begin{remark}
It is proved in~\cite{GLWZ06} that the conservativity problem is \textsc{NexpTime}-complete in the modal systems $\mathsf{K, S5}$ and that it is in $\textsc{ExpSpace}$ and \textsc{NexpTime}-hard in $\mathsf{S4}$. All these logics have the interpolation property~\cite{Mak1}, so according to the results of this section, the same lower and upper bounds apply to the admissibility problem for $\Pi_2$-rules.
\end{remark}

\section{Uniform Interpolants}\label{sec:globalUI}

We now present another strategy to determine the admissibility of $\Pi_2$-rules via (local and global) \emph{uniform interpolation}: this is a strengthening of ordinary interpolation.

\begin{definition}
A \emph{uniform local pre-interpolant} of a formula $\phi(\ux, \uy)$ wrt the variables $\ux$ is a formula  
$\exists_{\ux}^l \phi$ such that: (i) in  $\exists_{\ux}^l \phi$ at most the variables $\uy$ occur;
(ii) for every formula $\psi(\uy, \uz)$, we have 
\begin{equation}\label{eq:UI}
 \vdash_{\cS} \exists_{\ux}^l \phi \to \psi~~{\rm iff}~~\vdash_{\cS} \phi \to \psi~~.
\end{equation}
\end{definition}

Since $\vdash_{\cS} \phi \to\exists_{\ux}^l \phi$ holds as a special case by taking $\psi$  equal to $\phi$, if a uniform local pre-interpolant exists for every $\phi$, then $\cS$ has the interpolation property. 
The converse implication holds in case the logic is locally tabular, because in that case one can take as $\exists_{\ux}^l \phi$ the conjunction of all formulas $\psi(\uy)$ which are implied by $\phi$.\footnote{
It should be clear however that local tabularity alone is not sufficient for existence of local uniform pre-interpolants, because the set of formulas of the kind $\psi(\uy, \uz)$ implied by $\phi$ is not finite, modulo equivalence in $\cS$:
ordinary interpolation is needed for  the conjunction over the finite set  $\{\psi(\uy)\mid \phi\to \psi\}$ to be a local uniform pre-interpolant.
}
Uniform interpolants rarely exist, they are well-studied in the literature \cite{GZ}.

In case uniform local pre-interpolants exist, we  have a trivial criterion for conservativity (and consequently for admissibility of $\Pi_2$-rules). This was already pointed out in~\cite{GLWZ06}, we repeat the same observation in our context.

\begin{proposition}\label{prop:ui}
 If the local uniform pre-interpolant $\exists_{\up}^l F$ exists, then a
$\Pi_2$-rule $\rho$ of the form
 \[
\inference{F(\underline{\varphi}/\ux,\up) \impl \chi}
   {G(\underline{\varphi}/\ux) \impl \chi}
\]
is admissible in $\cS$ iff 
 \begin{equation}\label{eq:localUI}
  \vdash_{\cS}G\to \exists_{\up}^l F.
 \end{equation}
\end{proposition}

\begin{proof}
 Combining Theorem~\ref{thm:cons} and~\eqref{eq:UI} above, it is sufficient to observe that
 $\vdash_{\cS}G\to \exists_{\up}^l F$ holds precisely iff  $G(\ux)\wedge F(\ux, \up)$ is a conservative extension  of $G(\ux)$ in $\cS$.
\end{proof}

In the remaining part of this section, we will be interested in \emph{global} uniform interpolants: these are obtained by replacing in~\eqref{eq:UI}, local consequence relation by global consequence relation. In more detail:
\begin{definition}
A \emph{uniform global pre-interpolant} of a formula $\phi(\ux, \uy)$ wrt the variables $\ux$ is a formula  
$\exists_{\ux}^g \phi$ such that: (i) in  $\exists_{\ux}^g \phi$ at most the variables $\uy$ occur;
(ii) for every formula $\psi(\uy, \uz)$, we have 
\begin{equation}\label{eq:GUI}
 \exists_{\ux}^g \phi \vdash_{\cS}  \psi~~{\rm iff}~~ \phi \vdash_{\cS} \psi.
\end{equation}
\end{definition}

Since we will exclusively be interested in global uniform interpolants, we will write $\exists_{\ux}\phi$ for $\exists_{\ux}^g \phi$ and when we  speak of uniform interpolants, we will always mean global uniform interpolants. 
We write $\forall_{\ux} \phi$ for 
$\neg \exists_{\ux}\neg  \phi$; notice that,  for every formula $\psi(\uy, \uz)$, we have
\begin{equation}\label{eq:UUI}
 \psi\vdash_{\cS}\forall_{\ux} \phi   ~~{\rm iff}~~\psi  \vdash_{\cS}\phi.
\end{equation}
Existence of local and global interpolants are independent: for instance, in the basic modal logic $\mathsf{K}$ local uniform interpolants exist and global uniform interpolants do  not exist \cite{KM19}. For the converse phenomena, it is sufficient to recall once more the counterexamples from~\cite{Mak,Mak1} (notice that, in the locally tabular case, uniform local/global interpolants collapse to the corresponding ordinary local/global interpolants). 

There is no evident reason why global uniform interpolants could lead to recognizability of $\Pi_2$-rules, as it happens for the local uniform interpolants case (see Proposition~\ref{prop:ui} above). We will nevertheless show that this is the case when global uniform interpolants are accompanied by \emph{universal modalities}.

In the remaining part of this section and in the next section, we assume that our modal system $\cS$ has a universal modality.
We may view our modal signature $\Sigma$ as a first-order signature and $\Sigma$-formulas as terms of this signature. For a modal system $\cS$, an $\cS$-algebra is a Boolean algebra with operators (one operator of suitable arity for each $\Box\in \Sigma$) satisfying 
$\phi=\top$ for every   $\cS$-axioms $\phi$.
We call an $\cS$-algebra \emph{simple} iff the universal first-order sentence 
$\forall x \,(\univ x=\tt \vee \univ x=\ff)$ holds. This agrees with the standard definition from universal algebra, because it can be shown that congruences in an $\cS$-algebra bijectively correspond to 
\emph{$\univ$-filters}, i.e.,  filters $F$ satisfying the additional condition that $a\in F$ implies $\univ a\in F$. 
A standard Lindenbaum construction proves the \emph{algebraic completeness theorem}, namely that for every $\phi$ we have $\cS\vdash \phi$ iff the identity $\phi=\tt$ holds in all $\cS$-algebras (and hence iff $\phi=\tt$ holds in all simple $\cS$-algebras, because $\cS$-algebras are a discriminator variety \cite{jips:disc93}).

We need also some definitions from unification theory. A formula $\univ\phi(\ux)$ is said to be \emph{unifiable} in $\cS$ iff there is a substitution $\sigma$ mapping the $\ux:=x_1, \dots, x_n$ to some formulas $\sigma(\ux):= \sigma(x_1), \dots, \sigma(x_n)$ such that $\vdash_{\cS} \univ\phi(\sigma(\ux)/\ux)$. Such a substitution is said to be a \emph{unifier} of $\phi$ in $\cS$. The unifier is said to be \emph{projective} iff we have in addition
\begin{equation*}
 \univ\phi(\ux)\vdash_{\cS} \bigwedge_i (x_i \leftrightarrow \sigma(x_i)).
\end{equation*}
Notice that, by Lemma~\ref{lem:repl}, this implies
\begin{equation}\label{eq:pu}
 \univ\phi(\ux)\vdash_{\cS} \psi \leftrightarrow \psi(\sigma(\ux)/\ux).
\end{equation}
for every formula $\psi(\ux)$.

\begin{theorem}\label{thm:prun}
 Let $\cS$ have a universal modality. Then every formula $\univ\phi$ which is unifiable in $\cS$ has a projective unifier.\footnote{
 Despite this strong result,
 unifiability itself turns out to be undecidable for common modal systems with a universal modality, see~\cite{unifundec}.
 }
\end{theorem}

\begin{proof}
This result follows from the proof of the unitarity of unification in
  discriminator varieties, see~\cite{Burris} or also~\cite{Ry}. 
  We give here a direct simple proof, obtained by a slight generalization of an argument from~\cite{usurvey}.
  
Let $\univ\phi$ be unifiable; then there is a substitution $\sigma$ such that  
$\vdash_{\cS} \univ\phi(\sigma(\ux)/\ux)$.
Consider the substitution $\pi$ mapping a variable $x$ to
$$
(\univ\phi\wedge x) \vee (\neg \univ \phi \wedge \sigma(x)).
$$
This substitution clearly enjoys the property~\eqref{eq:pu}, so we only need to check that it is a unifier for $\univ\phi$. Consider now a simple $\cS$-algebra $\cA$ and a valuation $V$ of the propositional formulas into the support of $\cA$. By induction, it is easy to see that for every formula $\psi$ we have the following:
\begin{itemize}
 \item[-] if $V(\univ \phi)=\top$, then $V(\pi(\psi))= V(\psi)$;
 \item[-] if $V(\univ \phi)=\bot$, then $V(\pi(\psi))= V(\sigma(\psi))$.
\end{itemize}
In particular, for $\psi=\univ\phi$, we have that if $V(\univ \phi)=\top$, then $V(\pi(\univ\phi))= V(\univ\phi)=\top$ and if $V(\univ \phi)=\bot$, then $V(\pi(\univ\phi))= V(\sigma(\univ\phi))=\top$. Thus, for every simple algebra $\cA$ and for every valuation $V$ on the support of $\cA$, we have that $V(\pi(\univ \phi))=\top$
(notice that in a simple $\cS$-algebra, the only elements of the kind $\univ a$ are just $\top$ and $\bot$). Since $\cS$ is a discriminator variety, it is generated by its simple algebras, hence we have 
that $\pi$ unifies $\univ \phi$. \end{proof}

We need a technical lemma, showing a `Beck-Chevalley' condition, namely that uniform interpolants are stable under substitution, in the following sense: suppose that the (global) uniform interpolant $\exists_{\ux}\phi$
of $\phi(\ux, \uy)$ exists. This is a formula in the variables $\uy:=y_1, \dots, y_m$, so that for a tuple of formulas  $\upsi:=\psi_1, \dots, \psi_m$, it makes  sense to consider the formula  $(\exists_{\ux}\phi)(\upsi/\uy)$. But then one can consider the formula 
$\phi(\ux,\upsi/\uy)$ and the uniform interpolant $\exists_{\ux}\phi(\ux,\upsi/\uy)$: if the $\upsi$ do not contain the $\ux$, the following lemma ensures that the two formulas are the same (modulo provable equivalence in $\cS$).

\begin{lemma}\label{lem:bc}
 Let the formulas $\upsi$ do not contain the variables $\ux$ and suppose that the uniform interpolant $\exists_{\ux}\phi$ of $\phi(\ux, \uy)$ exists. Then $(\exists_{\ux}\phi)(\upsi/\uy)$ is the uniform interpolant $\exists_{\ux}\phi(\ux,\upsi/\uy)$.
\end{lemma}

\begin{proof}
 We need to show that for every formula $\theta$ not involving the $\ux$, we have 
 $$
 (\exists_{\ux} \phi)(\upsi/\uy)\vdash_{\cS} \theta \quad {\rm iff} \quad \phi(\ux, \upsi/\uy) \vdash_{\cS} \theta.
 $$
Notice that since the $\uy$ are used just as placeholders for substitutions,
after a suitable renaming, 
we can freely suppose that the $\uy$ do not occur in $\theta$ and  in the $\upsi$. 

The left-to-right side comes from
$\phi(\ux, \uy)\vdash_{\cS} \exists_{\ux} \phi$, by 
applying the substitution mapping the $\uy$ to the $\upsi$.
For the other side, assume  that 
$\phi(\ux, \upsi/\uy) \vdash_{\cS} \theta$; by Lemma~\ref{lem:repl}, we have (supposing that $\upsi:=\psi_1, \dots, \psi_m$) 
$$
\bigwedge_i (\psi_i\leftrightarrow y_i) \vdash_{\cS}  \phi(\ux, \uy) \leftrightarrow  \phi(\ux,\upsi/\uy)
$$
hence also 
$$
\phi\wedge \bigwedge_i (\psi_i\leftrightarrow y_i) 
\vdash_{\cS} \theta
$$
and 
$$
\phi\vdash_{\cS}
M \left( \bigwedge_i (\psi_i\leftrightarrow y_i) \right) \to
 \theta
$$
(for a suitable modality $M$).
By the definition of a global uniform interpolant, we get 
$$
\exists_{\ux}\phi\vdash_{\cS}
M \left( \bigwedge_i (\psi_i\leftrightarrow y_i) \right) \to
 \theta
$$
and finally 
$$
(\exists_{\ux}\phi)(\upsi/\uy)\vdash_{\cS}
 \theta
$$
applying the substitution mapping the $\uy$ to the $\upsi$.
\end{proof}

Putting together Theorem~\ref{thm:prun} and Lemma~\ref{lem:bc}, we can prove our characterization of $\Pi_2$-rules in presence of a universal modality and of uniform global interpolants.

\begin{theorem}\label{thm:pi2+UI+UM}
 Suppose that $\cS$ has uniform global interpolants and a universal modality. Then the rule
 $\rho$ is admissible in $\cS$ iff 
 \begin{equation}\label{eq:admpi2}
  \vdash_{\cS} \univ \forall_{\up} ( F(\ux, \up) \to z) \to (G(\ux) \to z).
 \end{equation}
\end{theorem}

\begin{proof}
 Suppose first that~\eqref{eq:admpi2} holds and assume that the antecedent $F(\uphi/\ux,\up)\to \chi$
 of the rule is provable in $\cS$. Then, after applying~\eqref{eq:UUI} and necessitation, we have 
 $$
 \vdash_{\cS} \univ \forall_{\up}  (F(\uphi/\ux, \up) \to \chi),
 $$
 if we apply the substitution mapping the $\ux$ to the $\uphi$ and $z$ to $\chi$ to~\eqref{eq:admpi2}
 (recall also Lemma~\ref{lem:bc}), we obtain
 $$
 \vdash_{\cS} \univ  \forall_{\up} ( F(\uphi/\ux, \up) \to \chi) \to (G(\uphi/\ux) \to \chi)
 $$
 yieldying $\vdash_{\cS} G(\uphi/\ux) \to \chi$ by modus ponens. This shows that $\rho$ is admissible.
 
 Conversely, suppose that $\rho$ is admissible. Consider the formula
 \begin{equation}\label{eq:ant}
  \univ\forall_{\up}  (F(\ux, \up) \to z).
 \end{equation}
 This is unifiable (a unifier is the substitution mapping $z$ to $\top$ and mapping the remaining variables to themselves),\footnote{Notice that this argument requires Lemma~\ref{lem:bc} too, applied to the formula $F(\ux, \up) \to z$.}   hence it has a projective unifier according to Theorem~~\ref{thm:prun}. Let us now make this projective unifier explicit. The variables occurring  in~\eqref{eq:ant} are the $\ux$ and $z$. So, suppose that the unifier maps the $\ux$ componentwise into certain formulas $\uphi$ and $z$ to a certain formula $\chi$.
 Thus we have, according to the definition of a unifier
 $$
 \vdash_{\cS} \univ\forall_{\up}  (F(\uphi/\ux, \up) \to \chi).
 $$
 Applying the reflexivity axiom for $\univ$ and~\eqref{eq:UUI}, we get
 $\vdash_{\cS} F(\uphi/\ux, \up) \to \chi$, i.e.,
 $$
 \vdash_{\cS} G(\uphi/\ux) \to \chi
 $$
 by the admissibility of $\rho$. Applying~\eqref{eq:pu}, we obtain
 $$
 \univ\forall_{\up}  (F(\ux, \up) \to z)\vdash_{\cS}
 (G(\ux) \to z) \leftrightarrow
 (G(\uphi/\ux) \to \chi)
 $$
 and also 
 $$
 \univ\forall_{\up}  (F(\ux, \up) \to z)\vdash_{\cS}
 (G(\ux) \to z)
 $$
 by modus ponens. This implies~\eqref{eq:admpi2} by Lemma~\eqref{lem:dedrepl}(i) and the transitivity axiom for $\univ$. 
 \end{proof}

 We summarize the results of this section into the following:
 
 \begin{theorem}\label{thm:main} Suppose that $\cS$ is decidable.
  $\Pi_2$ rules are effectively recognizable in $\cS$ in case that either 
  \begin{description}
   \item[{\rm (i)}] $\cS$ has computable local uniform interpolants or
   \item[{\rm (ii)}] $\cS$ has a universal modality and computable global uniform interpolants.
  \end{description}
 \end{theorem}

The two conditions above are independent: in fact (i) applies also to modal systems (like $\mathsf{K}$) without the universal modality. On the one hand, in presence of a universal modality, the existence of uniform global interpolants is weaker than the existence of local uniform intepolants (if we have both local uniform interpolants and a universal modality, we can define $\exists^g_{\ux}
 \phi$ as $\exists^l_{\ux} \univ \phi$). On the other hand, it is easy to check that in case we have both local uniform interpolants and a universal modality, the conditions for admissibility supplied by Theorem~\ref{thm:pi2+UI+UM} and Proposition~\ref{prop:ui} are equivalent.  For verifying this notice that if ~\eqref{eq:admpi2} is provable, then \eqref{eq:localUI}
 is also provable, taking $\exists^l_{\ux} F$ as $z$ (because $\vdash_{\cS}  F\to \exists^l_{\ux} F$). The converse implication can be proved in `natural deduction style' as follows: assume $\univ \forall_{\up}^g ( F(\ux, \up) \to z)$ and $G(\ux)$: then, it is possible to deduce $F(\ux, \up) \to z$ by $\forall_{\up}^g ( F(\ux, \up) \to z)\vdash_{\cS} F(\ux, \up) \to z$ and Lemma~\ref{lem:dedrepl}(i); by existential quantifier introduction, you can get $\exists^l_{\ux} F(\ux, \up) \to z$ 
 and finally $z$ by \eqref{eq:localUI} and implication elimination (modus ponens).

\section{Model Completions}\label{sec:modelcompl}

Uniform global interpolants are closely connected to model completions~\cite{GZ,MetReg}. This connection paves an alternative way for recognizing admissibility of $\Pi_2$-rules via algebraic and semantic methods. 

\subsection{An admissibility criterion}

Before continuing, we need to recall a few  results (restated as Theorems~\ref{thm:completeness} ad~\ref{thm:criterion} below) from~\cite{BBSV}. Since, in order to adapt them to our context, we need a slight generalization of these results, we provide the proofs in full detail.

With each $\Pi_2$-rule $\rho$ 
(see Definition~\ref{defirules}),
we associate the following  $\forall\exists$-statement in the \emph{first-order} language of $\cS$-algebras:
\[
\Pi(\rho) := \forall \ux,z  \Big( G(\ux) \nleq z 
\Rightarrow
\exists \uy : F(\ux,\uy) \nleq z \Big).
\]

We call $T_\cS$ the equational first-order theory of \emph{simple non degenerate $\cS$-algebras} (an $\cS$-algebra is non degenerate iff $\bot \neq \top$).

\begin{theorem}\label{completeness_rules}\label{thm:completeness}
 For each set  $\Theta$ of $\Pi_2$-rules and each formula $\psi$, we have 
 $$T_\cS \cup \{\Pi(\rho)\mid \rho \in \Theta\} \models \psi=\top ~\Longleftrightarrow~ \vdash_{\cS+\Theta}\psi.$$
\end{theorem}

\begin{proof}
The right-to-left direction is a trivial induction on the length of a proof witnessing  $\vdash_{\cS+\Theta}\psi$. For the other direction, we need a modified version of Lindenbaum's construction. Suppose that $\not\vdash_{\cS+\Theta}\psi$.
For each rule $\rho_i\in \Theta$, we add a countably infinite set of fresh propositional
letters to the set of existing propositional letters. Then we build the Lindenbaum algebra $\cB$ over the expanded set of
propositional letters, where the elements are the equivalence classes $[\varphi]$ under provable equivalence in $\cS+\Theta$. Next we 
construct a maximal $\univ$-filter $M$ of $\cB$ such that
$ [\neg \univ \psi]\in M$ and for every rule $\rho_{i}\in \Theta$
\[
(\rho_i) \quad \inference{F_i(\underline{\varphi}/\ux,\up) \impl \chi}
   {G_i(\underline{\varphi}/\ux) \impl \chi}
\]
and formulas
$\underline{\varphi}, \chi$: 

\vskip 2mm\noindent
$(\dag)$~  if $[G_{i}(\underline{\varphi}) \to \chi] \not\in M$, then there is a tuple $\up$ such that
$[F_{i}(\underline{\varphi},\up) \to \chi] \not\in M$.

\vskip 2mm\noindent
To construct $M$, let $\Delta_0 :=  \{ [\neg\univ\varphi] \}$, a consistent set. We enumerate all formulas $\varphi$ 
as $(\varphi_k : k \in \mathbb{N})$ 
and all tuples $(i,\underline{\varphi},\chi)$ 
where $i \in \mathbb{N}$ and $\underline{\varphi},\chi$ are as in the particular rule $\rho_i$, 
and we build the sets $\Delta_0 \subseteq \Delta_1\subseteq \cdots \subseteq \Delta_n \subseteq \cdots $ as follows (notice that, according to the construction below,  for all $n$ and  $\theta\in \Delta_n$, we have $\vdash_{\cS+\Theta} \theta\leftrightarrow \univ \theta$).
\begin{itemize}
\item For $n=2k$, if $\not\vdash_{\cS+\Theta}\bigwedge \Delta_n \to \univ \varphi_k$, let 
$\Delta_{n+1}=\Delta_n \cup \{ \neg\univ \varphi_k \}$; otherwise let
$\Delta_{n+1} = \Delta_n$.
\item For $n=2k + 1$, let $(l,\underline{\varphi},\chi)$ be the $k$-th tuple. If 
$\not\vdash_{\cS+\Theta}\bigwedge \Delta_n \to  (G_l(\underline{\varphi}) \to \chi)$,
let $\Delta_{n+1}=\Delta_n \cup \{\neg \univ (F_l(\underline{\varphi},\up) \to \chi)\}$,
where $\up$ is a tuple of propositional letters for $\rho_l$ not occurring in $\underline{\varphi},\chi$, and any of $\theta$ with 
$\theta\in \Delta_n$ (we can take $\up$ from the countably infinite additional propositional letters which we have reserved for the rule 
$\rho_l$). Otherwise, let $\Delta_{n+1} = \Delta_n$.
\end{itemize}
Let $M$ be 
$$
\{ ~[\theta] \mid ~\text{there are}~ \theta_1, \dots, \theta_n\in \bigcup_{n \in \mathbb{N}} \Delta_i~ \text{such that}~ \vdash_{\cS+\Theta}\theta_1\wedge\cdots \wedge \theta_n \to \theta \}.
$$
That $M$ is a proper $\univ$-filter not containing $[\psi]$ follows from the fact that $\not\vdash_{\cS+\Theta}\bigwedge \Delta_n \to \bot$. This is clear for $n=0$ and for any positive even $n$. For odd $n=2k+1$, suppose that $\vdash_{\cS+\Theta}\bigwedge \Delta_k  \to \univ (F_l(\underline{\varphi},\up) \to \chi) $ and that 
$\not\vdash_{\cS+\Theta}\bigwedge \Delta_k \to  (G_l(\underline{\varphi}) \to \chi)$. Then, by the reflexivity  axiom $\univ \phi \to \phi$, we have $\vdash_{\cS+\Theta} F_l(\underline{\varphi},\up) \to(\bigwedge\Delta_k    \to \chi)$
and also 
(applying the rule $\rho_l$ of the $k$-th tuple)
$\vdash_{\cS+\Theta} G_l(\underline{\varphi}) \to(\bigwedge\Delta_k    \to \chi)$, yielding a contradiction.

Also, by the even steps of 
the construction of the sets $\Delta_n$, it contains either $[\univ \theta]$ or $[\neg\univ \theta]$
for every $\theta$, thus $M$ is a maximal $\univ$-filter.
Finally, the odd steps of the construction of the sets $\Delta_n$ ensure that $M$ satisfies (\dag):
in fact, if $[G_{i}(\underline{\varphi}) \to \chi] \not\in M$, then by step $n=2k + 1$, we have $[\neg \univ (F_l(\underline{\varphi},\up) \to \chi)]\in M$ and if also $[F_{i}(\underline{\varphi},\up) \to \chi]\in M$, then $[\univ(F_{i}(\underline{\varphi},\up) \to \chi)]\in M$ (because $M$ is a $\univ$-filter) and so $M$ would not be proper, a contradiction.
Therefore, we can conclude that $M$ satisfies
all the desired properties.

By (\dag), the quotient of $\cB$ by $M$ satisfies each $\Pi(\rho_i)$. This quotient is a simple algebra, because $M$ is maximal as a $\univ$-filter.  Moreover, since $[\neg\univ \psi]\in M$, we have that $[\neg\univ \psi]$
maps to $\top$, so $[\univ \psi]$ maps to $\bot$ in the quotient. Thus, $[\varphi]$ does not map to $\top$ in the quotient,
and hence 
$T_\cS \cup \{\Pi(\rho)\mid \rho \in \Theta\}\not \models \psi=\top$.
\end{proof}

We will use $\&$, $\mathtt{or}$, $\foneg$, and $\Rightarrow$ to denote first-order connectives in order to distinguish them from Boolean algebra operations. $\bigwedge$ and $\bigvee$ will denote finite first-order conjunctions and disjunctions.

\begin{definition} \label{theencoding}
Given a quantifier-free first-order formula $\Phi(\ux)$
in the language of $\cS$-algebras, 
we associate with
it the term (aka the propositional modal formula)
$\Phi^\ast(\ux)$ as follows:
\begin{align*}
(t(\ux) = u(\ux))^\ast &= \univ (t(\ux)\biimpl u(\ux)) \\
(\foneg\Psi)^\ast(\ux) &= \neg \Psi^\ast(\ux) \\
(\Psi_1(\ux) \fowedge \Psi_2(\ux))^\ast &= \Psi_1^\ast(\ux) \cong
\Psi_2^\ast(\ux).
\end{align*}
\end{definition}

The following lemma is immediate:

\begin{lemma} \label{belleformule}
Let $\cB$ be a simple $\cS$-algebra and let 
 $\Phi(\ux)$ be a quantifier-free formula. Then 
 we have $$\cB\models \Phi(\ua/\ux)\ \mbox{iff}\  \cB\models (\Phi(\ua/\ux))^\ast=\top,$$ 
 for every  tuple $\ua$ from $\cB$.
\end{lemma}

\begin{theorem}
[Admissibility Criterion] \label{thm:criterion}
A $\Pi_2$-rule $\rho$ is admissible in $\cS$ iff for each 
simple $\cS$-algebra $\cB$ there is a simple $\cS$-algebra $\cC$
such that $\cB$ is a substructure  of $\cC$ and $\cC \models \Pi(\rho)$.
\end{theorem}

\begin{proof}
 ($\Rightarrow$) Suppose that the rule $\rho$ 
 \[
(\rho) \quad \inference{F(\underline{\varphi}/\ux,\up) \impl \chi}
   {G(\underline{\varphi}/\ux) \impl \chi}
\]
 is admissible in $\cS$. It is sufficient to show that there exists a model $\cC$ of the theory
\[
T = T_{\cS} \cup \{ \Pi(\rho) \} \cup \Delta(\cB)
\]
where $\Delta(\cB)$ is the diagram of $\cB$ \cite[p.~68]{CK}.
Suppose for a contradiction that $T$ has no models, hence is inconsistent. Then, by compactness, there exists a quantifier-free
first-order formula $\Psi(\ux) $ and a tuple $\ux$ of variables corresponding to some $\ua \in \cB$
such that
\[
T_{\cS} \cup \{ \Pi(\rho) \} \models \foneg \Psi(\ua/\ux) \text{ and } \cB \models \Psi(\ua/\ux).
\]
By Theorem~\ref{completeness_rules},
$\cS + \rho$ is complete with respect to the simple $\cS$-algebras  satisfying $\Pi(\rho)$.
Therefore, by Lemma~\ref{belleformule}, we have $T_{\cS} \cup \{ \Pi(\rho) \} \models (\foneg \Psi(\ux))^\ast=\top$ and also
$\vdash_{\cS + \rho} (\foneg \Psi(\ux))^{*}$, where $(-)^{*}$ is the translation given in Definition~\ref{theencoding}.
By admissibility,
$\vdash_{\cS} (\foneg \Psi(\ux))^{*}$. 
Thus, for the valuation $v$ 
into $\cB$ 
that maps $\ux$ to $\ua$, we have $v((\foneg\Psi(\ux))^{*}) = 1$, so $v((\Psi(\ux))^{*}) = 0$. This contradicts the 
fact that $\cB \models \Psi(\ua/\ux)$. Consequently, $T$ must be consistent, and hence it has a model.

($\Leftarrow$) Suppose $\vdash_{\cS} F(\underline{\varphi},\up)\to \chi$ with $\up$ not occurring in
$\underline{\varphi},\chi$. Let $\cB$ be a simple $\cS$-algebra  and let $v$ be a valuation on $\cB$. By assumption, there is a simple $\cS$-algebra $\cC$
such that $\cB$ is a substructure of $\cC$ and $\cC \models \Pi(\rho)$. Let
$i : \cB \hookrightarrow \cC$ be the inclusion. Then $v':= i \circ v$ is a valuation on $\cC$.
For any $\uc \in C$,
let $v''$ be the valuation that agrees with $v'$ except for the fact that it maps the $\up$ into the $\uc$.
Since $\vdash_{\cS} F(\underline{\varphi}/\ux,\up) \to \chi$, by the 
algebraic completeness theorem\footnote{This is Theorem~\ref{completeness_rules} for $\Theta=\emptyset$.}
 we have
$v''(F(\underline{\varphi}/\ux,\up) \to \chi) = \top$. This means that for all $\uc \in \cC$, 
we have $F(v'(\underline{\varphi}),\uc) \leq v'(\chi)$. Therefore, $\cC\models \forall \uy
\Big( F(v'(\underline{\varphi}),\uy) \leq v'(\chi) \Big)$. Since $\cC \models \Pi(\rho)$, we have
$\cC \models G(v'(\underline{\varphi})) \leq v'(\chi)$. Thus, as $G(v'(\underline{\varphi})) \leq v'(\chi)$ holds in $\cC$,
we have that $G(v(\underline{\varphi})) \leq v(\chi)$ holds in $\cB$. Consequently, $v(G(\underline{\varphi}) \to \chi) = \top$.
Applying the algebraic completeness theorem again yields that $\vdash_{\cS} G(\underline{\varphi}) \to \chi$ 
because $\cB$ is arbitrary, 
and hence
$\rho$ is admissible in $\cS$.
\end{proof}

\subsection{Admissibility and Model Completeness}

We now investigate the connections between admissibility and model completions.

\begin{theorem}\label{thm:mainadm} Suppose that the universal theory $T_{\cS}$ has a model completion $T_{\cS}^\star$. Then 
 a $\Pi_2$-rule $\rho$ is admissible in  $\cS$ iff $T_{\cS}^\star\models \Pi(\rho)$.
\end{theorem}

\begin{proof} 
 Applying Theorem~\ref{thm:criterion},
 we show that
 $T_{\cS}^\star\models \Pi(\rho)$ holds iff  every simple $\cS$-algebra $\cB$  can be embedded into some simple $\cS$-algebra $\cC$ which  satisfies $\Pi(\rho)$.  This is shown below using the fact that $\Pi(\rho)$ is a $\Pi_2$-sentence. Recall that models of $T_{\cS}^\star$  are just the existentially closed simple $\cS$-algebras  (see~\cite[Proposition 3.5.15]{CK}).
 
 Suppose for the left to right direction that  $T_{\cS}^\star\models \Pi(\rho)$ holds and let $\cB$ be any simple $\cS$-algebra. Then $\cB$ embeds into an existentially closed simple $\cS$-algebra $\cC$ (this is a general model-theoretic fact~\cite{CK}). 
Thus, $\cC$ is a model of $T_{\cS}^\star$ and hence $\cC \models \Pi(\rho)$.
 
Conversely, suppose that every simple $\cS$-algebra $\cB$  can be embedded into some simple $\cS$-algebra $\cC$ which  satisfies $\Pi(\rho)$. Pick $\cB$ such that $\cB\models T_{\cS}^\star$  and let $\Pi(\rho)$ be $\forall \ux \exists \uy H(\ux,\uy)$, where $H$ is quantifier free. Let $\ub$ be a tuple from the support of $\cB$.
Let $\cC$ be an extension of $\cB$ such that $\cC \models \Pi(\rho)$. Then $\cC \models \exists \uy H(\ub,\uy)$.
 As $\cB$ is existentially closed, this immediately entails that $\cB\models \exists \uy H(\ub,\uy)$.
 Since the $\ub$ was arbitrary, we conclude that $\cB\models \Pi(\rho)$, as required.
\end{proof}

\begin{remark}
Theorems~\ref{thm:pi2+UI+UM} and~\ref{thm:mainadm} are in fact equivalent statements: indeed the existence of global uniform interpolants and the existence of a model completion for $T_{\cS}$ are equivalent statements (as it can be deduced from slight modifications of the results in~\cite{GZ,MetReg}) and if one considers how quantifiers are eliminated in $T_{\cS}^\star$ via global uniform interpolants~\cite{GZ}, one can 
translate the statements of Theorems~\ref{thm:pi2+UI+UM} and~\ref{thm:mainadm} into each other. 
We nevertheless point out that the two theorems are proved via completely different tools (namely unification theory and model theoretic techniques): 
this is  quite a notable fact. 
\end{remark}

According to Theorem~\ref{thm:mainadm}, checking whether a $\Pi_2$-rule is admissible now amounts to checking whether  $T_{\cS}^\star\models \Pi(\rho)$ holds. The latter can be done via quantifier elimination in $T_{\cS}^\star$. We give sufficient conditions for this to be effective. 

\begin{corollary}\label{maincoro} Let $\cS$ be decidable and locally tabular. Assume also that  simple $\cS$-algebras enjoy the amalgamation property.
 Then admissibility of $\Pi_2$-rules in $\cS$
 is effective.
\end{corollary}

\begin{proof}
Local tabularity of $\cS$ implies local finiteness\footnote{Recall that a class of algebras is \emph{locally finite} if every finitely generated algebra in this class if finite, see \cite[Section 14.2]{CZ97} for the connection between local finiteness and local tabularity.}  of $T_{\cS}$. For universal locally finite theories in a finite language, amalgamability is a necessary and sufficient condition for existence of a model completion~\cite{LIP,wheeler}. 
Quantifier elimination in  $T_{\cS}^\star$ is effective  because there are only finitely many non-equivalent formulas in a fixed finite number of variables, because of Lemma~\ref{belleformule}  and because of the following folklore lemma.
\end{proof}

\begin{lemma}\label{lem:qe}
 The quantifier-free formula $R(\ux)$ provably equivalent in $T_{\cS}^\star$ to an existential formula $\exists \uy H(\ux,\uy)$ is the strongest quantifier free formula $G(\ux)$ implied 
 (modulo $T_{\cS}$) by $H(\ux,\uy)$.
\end{lemma}

\begin{proof}
 Recall that $T_{\cS}$ and $T_{\cS}^\star$ are co-theories~\cite{CK}, i.e.~they prove the same universal formulas. Thus we have the following chain of equivalences:
 $$
 \begin{array}{c}
  T_{\cS} \vdash H(\ux, \uy) \to G(\ux)
  \\ 
  \hline 
  T_{\cS}^\star \vdash H(\ux, \uy) \to G(\ux)
  \\ 
  \hline 
  T_{\cS}^\star \vdash \exists \uy H(\ux, \uy) \to G(\ux)
  \\ 
  \hline 
  T_{\cS}^\star \vdash R(\ux) \to G(\ux)
  \\ 
  \hline 
  T_{\cS} \vdash R(\ux) \to G(\ux)
  \end{array}
  $$
 yielding the claim.
\end{proof}

The usefulness of Corollary~\ref{maincoro} lies in the fact  that its only real requirement is the amalgamation property, besides local tabularity. Whenever local tabularity holds, finitely presented algebras are finite, thus it is sufficient to establish amalgamability for \emph{finite} algebras:
in fact, two algebras $\cB_1, \cB_2$ amalgamate over a common subalgebra $\cA$ iff
$T_{\cS}\cup \Delta(\cB_1)\cup \Delta(\cB_2)$ is consistent iff (by compactness and local finiteness) there are amalgamating finite subalgebras $\cB^0_1, \cB^0_2$ 
of $\cB_1, \cB_2$, respectively.
Whenever a ``useful'' duality is established, amalgamation of finite algebras turns out to be equivalent to 
co-amalgamation
for finite frames, which is usually much easier to check.
We will now give a few simple examples and counterexamples.

\begin{example}
 If the modal signature contains only the universal modality $\univ$, we have the locally tabular  logic $\mathsf{S5}$. Finite simple non degenerate $\mathsf{S5}$-algebras are dual to finite nonempty sets and onto maps, for which 
co-amalgamation
 trivially holds (by standard pullback construction), see, e.g., \cite[Thm.~14.23]{CZ97}.
\end{example}

\begin{example}
 The logic of the difference modality~\cite{derijke,Ven93} has in addition to the global modality a unary operator $D$ subject to the axioms
 $$
 \univ \phi \biimpl (\phi \cong \ne D\ne \phi),\qquad \phi \impl \ne D \ne D\phi, \qquad DD\phi \impl \phi \dis D \phi. 
 $$
 This logic axiomatizes Kripke frames where the accessibility relation is inequality.
 Local finiteness can be established for instance by the method of irreducible models~\cite{Girr}. Amalgamation however fails. To see this,  notice that the simple frames for this logic are sets endowed with a 
 relation $E$ such that 
 $w_1\neq w_2\to w_1Ew_2$.
Now let $X = \{x_1,\dots, x_5\}$, $Y = \{y_1,\dots, y_5\}$ and $Z = \{z_1, z_2\}$. Let $x_i E_X x_j$ iff $i\neq j$ for $1\leq i, j\leq 5$, $y_i E_Y y_j$ iff $i\neq j$ for $1\leq i, j\leq 5$
and $z_i E_Z z_j$ for $i, j =1, 2$. Let also $f: X\to Z$ and $g: Y\to Z$ be such that $f(x_1) = f(x_2)= f(x_3) = g(y_1)= g(y_2) = z_1$ and  $f(x_4) = f(x_5)= g(y_3) = g(y_4)= g(y_5) = z_2$.
Then it is easy to see that $f$ and $g$ are p-morphisms. If a 
co-amalgam
exists, then there must exist a frame $(U, E_U)$ and onto p-morphims $h: U\to X$ and $j:U\to Y$ such that 
$f\circ h = g\circ j$. However, an easy argument shows that  $U$ should contain more than $5$ points. Moreover, for $u,v\in U$ with $u\neq v$ we should have $uE_Uv$.  
But then there will be distinct points in $U$ mapped by $h$ to some $x_i$, which would entail that $x_i$ is reflexive, which is a contradiction. 
\end{example}

\section{Symmetric Strict Implication and Contact Algebras}\label{sec:ca}

In this section we apply the results of Section~\ref{sec:modelcompl} (in particular, Corollary~\ref{maincoro}) to the case of contact algebras. We
first review some material from~\cite{BBSV}. Let us consider the modal signature comprising, besides the universal modality $\univ$, a binary operator $\simpa$, which we call \emph{strict implication}, subject to the following axioms (we keep the same enumeration as in~\cite{BBSV} and add axiom (A0) which is seen as a definition of $\univ$ in~\cite{BBSV}). 

\begin{ourlist}
    \item[(A0)] $\univ \varphi \biimpl (\top \simpa \varphi)$,
    \item[(A1)] $(\bot \simpa \varphi) \cong (\varphi \simpa \top)$,
    \item[(A2)] $[(\varphi \dis \psi) \simpa \chi] \biimpl [(\varphi \simpa \chi) \cong (\psi \simpa \chi)]$,
    \item[(A3)] $[\varphi \simpa (\psi \cong \chi)] \biimpl [(\varphi \simpa \psi) \cong (\varphi \simpa \chi)]$,
    \item[(A4)] $(\varphi \simpa \psi) \impl (\varphi \impl \psi)$,
    \item[(A5)] $(\varphi \simpa \psi) \biimpl (\ne \psi \simpa \ne \varphi)$,
    \item[(A8)] $\univ \varphi \impl \univ \univ \varphi$,
    \item[(A9)] $\ne\univ \varphi \impl \univ\neg\univ \varphi$,
    \item[(A10)] $(\varphi \simpa \psi) \biimpl \univ (\varphi \simpa \psi)$,
    \item[(A11)] $\univ \varphi \impl (\neg \univ \varphi \simpa \bot)$.
\end{ourlist}
Inference rules are modus ponens (for $\to$) and necessitation (for $[\forall]$). It can be shown (see~\cite{BBSV}) that this system (called \emph{symmetric strict implication calculus} $\sf S^2IC$) matches our requirements from Section~\ref{sec:preliminaries}.\footnote{ 
Strictly speaking, since $\simpa$ turns disjunctions into conjunctions in the first argument, to match those requirements we should replace the connective $\simpa$ with an equivalent binary modality $\Box$ related to $\simpa$ via the definition $x\simpa y := \Box[\neg x, y]$.
} 

We recall  that a \emph{symmetric strict implication algebra} ($\sf S^2I$-algebra for short) is a pair $\mathcal{B} = (B, \simpa)$, where $B$ is a Boolean algebra and $\simpa : B \times B \to B$ a binary operation validating the axioms 
(A0)-(A11)  \cite[Section 3]{BBSV}. Then axioms (A0), (A8)-(A11)  yield  that $[\forall]: B\to B$ is an $\sf S5$-operator on $B$ such that for each $a\in B$ we have $[\forall] a = 1 \simpa a$. Then the variety of $\sf S^2I$-algebras is semi-simple (every subdirectly irreducible algebra is simple) and simple $\sf S^2I$-algebras are those 
$\sf S^2I$-algebras $\mathcal{B} = (B, \simpa)$  where we have that $a\simpa b$ is either $0$ or $1$. This entails that $\sf S^2IC$ is locally tabular (in algebraic terms, the variety of $\sf S^2I$-algebras is locally finite). 
For the proofs of all these facts we refer to \cite[Section 3]{BBSV}.
Thus, in a simple non-degenerate $\sf S^2I$-algebra, the operation $\simpa$ is in fact the characteristic function of a binary relation. Given a simple $\sf S^2I$-algebra $\cB$ we define $\prec$ by setting
\[
a \prec b \mbox{ iff } a \simpa b = 1.
\]
Then $\prec$ satisfies the following axioms:
\begin{ourlist}
    \item[(S1)] $0 \prec 0$ and $1 \prec 1$;
    \item[(S2)] $a \prec b,c$ implies $a \prec b \land c$;
    \item[(S3)] $a,b \prec c$ implies $a \lor b \prec c$;
    \item[(S4)] $a \leq b \prec c \leq d$ implies $a \prec d$;
    \item[(S5)] $a \prec b$ implies $a \leq b$;
    \item[(S6)] $a \prec b$ implies $\neg b \prec \neg a$.
\end{ourlist}

Conversely, if $\prec$ is a binary relation on $B$ satisfying (S1)--(S6), we define $\simpa: B \times B \to B$ by
\[
a \simpa b =
\begin{cases}
    1 \ \mbox{ if } a \prec b \\
    0 \ \mbox{ otherwise. }
\end{cases}
\]
Then $(B, \simpa)$ is a simple $\sf S^2I$-algebra (i.e., satisfies (A0)--(A11) and $\simpa$ has values in $\{0,1\}$). Moreover,   
\[
[\forall] a =
\begin{cases}
    1 \mbox{ if } a = 1 \\
    0 \mbox{ if } a \neq 1.
    \end{cases}
\]

Finally, we note that this correspondence is one-to-one  \cite[Section 3]{BBSV}.

Non-degenerate Boolean algebras endowed with a relation $\prec$ satisfying the above conditions (S1)-(S6) are called \emph{contact algebras}.\footnote{It is more common to use in contact algebras the \emph{contact relation} $\delta$ \cite{Vak07}, which is given by $a\delta b$ iff $a\not\prec \neg b$. However, we stick with our notation to stay close to our main reference \cite{BBSV}.} The class of all contact algebras 
and the corresponding first-order theory are both
  denoted by $\sf Con$. 
The above considerations suggest  translations from the theory of simple $\sf S^2I$-algebras into the theory of contact algebras, and vice versa. We are interested in detailing the translations at the level of quantifier-free formulas. 

\begin{itemize}
 \item \emph{Translation $\tau_1$ from contact algebras to simple $\sf S^2I$-algebras}.  We let $\tau_1(t\prec u)$ to be $t \simpa u = 1$ and $\tau_1(t =u)$ to be $t=u$; the translation $\tau_1$ operates identically on Boolean connectives. 
 \item \emph{Translation $\tau_2$ from simple $\sf S^2I$-algebras to contact algebras}. We translate a quantifier-free formula $\phi$ in three steps. In the first step we eliminate all $\forall$ symbols using axiom (A0) and then we \emph{flatten} $\phi$  by repeatedly applying the following transformation:
 $$
 \phi \longmapsto \exists x~(x=t \simpa u \fowedge \phi(x/t \simpa u)).
 $$
 After this step and after moving the existential quantifiers to the front, the formula to be translated has the form 
 $$
 \exists x_1, \ldots, x_n ~\left(\bigwedge_{i=1}^n (x_i=t_i\simpa u_i)\fowedge \psi\right)
 $$
 where $\psi$ is a quantifier-free formula in the language of Boolean algebras.
 In the second step, we translate this formula into the following formula in the language of contact algebras
 \begin{equation}\label{eq:firsttranslation}
  \exists x_1, \ldots, x_n ~\left(\bigwedge_{i=1}^n [(x_i=1 \fowedge  u_i\prec t_i)
  \fovee (x_i=0 \fowedge u_i\not\prec t_i)] 
  \fowedge \psi\right).
 \end{equation}
 In the last step, we apply distributivity law to~\eqref{eq:firsttranslation}, thus obtaining an exponentially large disjunction; from each disjunct, the existential quantifiers can be removed by 
replacing $x_i$ with $1$ or $0$.
The final result will be our $\tau_2(\phi)$.
\end{itemize}

The following proposition follows from the above considerations:
\begin{proposition}\label{prop:translation} Let $T_{S^2I}$ be the theory of simple $\sf S^2I$-algebras.
 For every quantifier-free formulas $\phi_1, \phi_2$ in the languages of contact and of simple $\sf S^2I$-algebras, respectively, we have that:
 \begin{description}
  \item[{\rm (i)}] ${\sf Con}\models \phi_1$ implies $T_{S^2I}\models \tau_1(\phi_1)$;
  \item[{\rm (ii)}] $ T_{S^2I} \models \phi_2$ implies ${\sf Con} \models \tau_2(\phi_2)$;
  \item[{\rm (iii)}] $ T_{S^2I} \models \phi_1\Leftrightarrow \tau_2(\tau_1(\phi_1))$ and  ${\sf Con} \models \phi_2\Leftrightarrow \tau_1(\tau_2(\phi_2))$.
 \end{description}
\end{proposition}

Since, as outlined above, 
the theory of non degenerate simple $\sf S^2I$-algebras is essentially the same (in fact, it is a syntactic variant) as the universal 
theory $\sf Con$ of contact algebras, we shall investigate the latter in order to apply Corollary~\ref{maincoro}.\footnote{
Notice also that computing quantifier elimination in the model completions  (once we proved that such model completions exist by Corollary~\ref{maincoro}) commutes with the translations, by Proposition~\ref{prop:translation} and Lemma~\ref{lem:qe}. This observation will be  used in Subsection~\ref{subsec:complexity}.
}
What we have to show in order to check the hypotheses of such a corollary is just that  $\sf Con$ is amalgamable.

To prove amalgamability, we need a duality theorem. In~\cite{BBSV16,Cel01,DV17} a duality theorem is established for the category  of contact algebras and \emph{$\prec$-maps} (a map $\mu: (\cB, \prec) \to (\cC, \prec)$ among contact algebras is said to be a $\prec$-map iff it is a Boolean homomorphism such that $a\prec b$ implies $\mu(a)\prec \mu(b)$). We will make use of that theorem but will modify it, because for amalgamation we need a duality for contact algebras and \emph{embeddings} in the model theoretic sense   
(this means that an embedding is an injective map that not only \emph{preserves} but also \emph{reflects} the relation $\prec$). We first recall the duality theorem of~\cite{BBSV16}, giving just minimal information that is indispensable for our purposes.

We say that a binary relation $R$ on a topological space $X$ is \emph{closed} if $R$ is a closed subset of $X\times X$ in the product
topology. Let $\mathsf{StR}$ be the category having (i) as objects the pairs $(X,R)$, where $X$ is a (non empty) Stone space and $R$ is a closed, reflexive and symmetric relation on $X$, and (ii) as arrows the continuous  maps $f:(X,R) \to (X',R')$ which are \emph{stable} (i.e.~such that $xRy$ implies $f(x)R' f(y)$ for all points $x,y$ in the domain of $f$). We define a contravariant functor
$$
(-)^\star: ~ \mathsf{StR}^{op} \to \mathsf{Con}_s
$$
into the category $\mathsf{Con}_s$ of contact algebras and $\prec$-maps as follows:
\begin{itemize}
 \item for an object $(X,R)$, the contact algebra $(X, R)^\star$ has $\mathsf{Clop}(X)$ the clopens of $X$ as carrier set (with union, intersection and complement as Boolean operations) and its relation $\prec$ is given by 
 $C\prec D$ iff $R[C]\subseteq D$ (here we used the abbreviation $R[C]=\{ x\in X\mid sRx~ {\rm for~some}~s\in C\}$);
 \item for 
 a stable 
 continuous map $f:(X,R) \to (X',R')$, the map $f^\star$ is the inverse image along $f$.
\end{itemize}

\begin{theorem}[\cite{BBSV16,DV17}] The functor $(-)^\star$ establishes an equivalence of categories.
\end{theorem}

We now intend to restrict this equivalence to the category $\mathsf{Con}_e$  of contact algebras and embeddings. To this aim we need to identify a suitable subcategory $\mathsf{StR}_e$ of $\mathsf{StR}$.
Now $\mathsf{StR}_e$ has the same objects as $\mathsf{StR}$, however a stable continuous map $f:(X_1,R_1) \to (X_2,R_2)$ is in $\mathsf{StR}_e$ iff  it satisfies the following additional condition:
\begin{equation}\label{eq:emb}
 \forall x, y \in X_2 \text{ } [x R_2 y \text{ } \Leftrightarrow \text{ } \exists \Tilde{x}, \Tilde{y} \in X_1 \text{ s.t. } f(\Tilde{x}) = x, \text{ } f(\Tilde{y}) = y \text{ } \fowedge \text{ } \Tilde{x} R_1 \Tilde{y}].
 \end{equation}
Notice that, since $R_2$ is reflexive, it turns out that a map satisfying~\eqref{eq:emb} must be surjective. We call the stable maps satisfying~\eqref{eq:emb} \emph{regular stable maps}, because it can be shown that these maps are just the regular epimorphisms in the category $\mathsf{StR}$.

\begin{theorem}\label{thm:duality} The functor $(-)^\star$, suitably restricted in its domain and codomain,  establishes an equivalence of categories between $\mathsf{StR}_e$ and  $\mathsf{Con}_e$.
\end{theorem}

\begin{proof}
 We need to show that $f$ satisfies condition~\eqref{eq:emb} above iff $f^{\star}$ is an embedding between contact algebras, i.e. iff it satisfies the condition 
 \begin{equation}\label{eq:emb1}
 (R_1 [f^{-1} (U)] \subseteq f^{-1} (V) \text{ } \Leftrightarrow \text{ } R_2 [U] \subseteq V) ~~\text{ } \forall \; U, V \in \mathsf{Clop}(X_2) 
 \end{equation}
where $\mathsf{Clop}(X_2)$ is the set of clopens of the Stone space $X_2$. We tranform condition~\eqref{eq:emb1} up to equivalence. First notice that, by the adjunction between direct and inverse image, \eqref{eq:emb1} is equivalent to 
\begin{equation}\label{eq:emb2}
 (f(R_1 [f^{-1} (U)]) \subseteq  V \text{ } \Leftrightarrow \text{ } R_2 [U] \subseteq V) ~~\text{ } \forall \; U, V \in \mathsf{Clop}(X_2).
 \end{equation}
 Now, in  compact Hausdorff spaces  closed relations and continuous functions map
closed sets to closed sets, hence  $f(R_1 [f^{-1} (U)])$ is closed and so, since clopens are a base for closed sets, \eqref{eq:emb2} turns out to be equivalent to
\begin{equation}\label{eq:emb3}
 (f(R_1 [f^{-1} (U)])   
 = 
 R_2 [U] ) ~~~~\text{ } \forall \; U\in \mathsf{Clop}(X_2).
 \end{equation}
 We now claim that~\eqref{eq:emb3}
 is equivalent to 
 \begin{equation}\label{eq:emb4}
  f(R_1 [f^{-1} (\{ x \})]) = R_2 [\{ x \}] ~~~~~~~~\text{ } \forall x\in X_2.
 \end{equation}
In fact,~\eqref{eq:emb4} implies ~\eqref{eq:emb3} because all operations $f(-), R[-], f^{-1}(-)$ preserve set-theoretic unions. The converse implication holds because of Esakia's lemma below applied to the down-directed  system $\{U\in \mathsf{Clop}(X_2) \mid x\in U\}$. Notice that Esakia's lemma applies because $f\circ R_1\circ f^{op}$ and $R_2$ are symmetric relations, since $R_1$ and $R_2$ are symmetric (here we view $f$ and $f^{-1}=f^{op}$ as relations via their graphs).

Now it is sufficient to observe that~\eqref{eq:emb4} is equivalent to the conjunction of~\eqref{eq:emb} and stability.
\end{proof}

We will now prove a version of Esakia's lemma for our spaces. Esakia's lemma normally speaks about the inverse of a relation $R$, but here we need a version which holds for $R$-images because our relation is symmetric. 

\begin{lemma}$($Esakia's lemma,~\cite[Lemma 3.3.12]{Esakia}$)$
Let $X$ be a compact Hausdorff space and $R$ a point-closed\footnote{A binary relation $R$ on a topological space $X$ is said to be \emph{point-closed} if $\forall x\in X$ $R[x]$ is closed in $X$.
A closed relation in  a compact Hausdorff space maps closed sets to closed sets via $R[-]$, hence it is point-closed.
} symmetric binary relation on $X$. Then for each downward directed
family $\mathcal{C} = $\{$ C_i $\}$_{i\in I} $ of nonempty closed subsets of $X$, we have $R [\underset{i\in I}{\bigcap}C_i] = \underset{i\in I}{\bigcap}R[C_i]$.
\end{lemma}
\begin{proof} The inclusion 
$R[\underset{i\in I}{\bigcap} C_i] \subseteq \underset{i\in I}{\bigcap} R[C_i]$ is trivial.
Now suppose $x\in \underset{i\in I}{\bigcap} R[C_i]$. Then $x\in R[C_i]$ for each $C_i$ and, by symmetry, $R[x]\cap C_i$  is nonempty for each $i\in I$. But as $C_i$'s are downward directed, all the finite intersections $R[x] \cap C_{i_1}\cap ...\cap C_{i_n}$ (with $i_j\in I$ for $j\in \{ 1, ..., n \}$) are nonempty. By compactness, the infinite intersection (which equals $R[x] \cap \underset{i\in I}{\bigcap} C_i$) is nonempty and so, by symmetry, $x\in R[\underset{i\in I}{\bigcap} C_i]$.
\end{proof}

Whenever there is a regular stable map $f: (Y, R') \longrightarrow (X,R)$, we say that $(Y,R')$ \emph{covers} $(X,R)$. The following lemma gives an interesting example of a cover and will be useful in Subsection~\ref{subsec:complexity} below.
Let us call \emph{contact frames} the objects of $\mathsf{StR}$. A \emph{singleton} in a contact frame $(Y,R)$ is a point $y\in Y$ such that 
$R[y]=\{ y \}$.
A contact frame  $(Y,R)$ is said to be a \emph{1-step contact frame} iff it does not contain singletons and it satisfies the following condition for all $x,y,z\in Y$:
\begin{equation}
 xRy \fowedge yRz ~~\Rightarrow~~(x=y \fovee y=z \fovee x=z).
\end{equation}
Thus the points in a 1-step contact frame can be partitioned into  2-element subsets
$\{y_1,y_2\}$ such that $y_1\neq y_2$ and the only elements accessible from $y_i$ ($i=1,2$) are $\{y_1,y_2\}$.

\begin{lemma}\label{lem:1step}
 Every finite contact frame $(X,R)$ is covered by a 1-step contact frame of at  most quadratic size.
\end{lemma}

\begin{proof} We first get rid of singletons by `duplicating' them: this means that we move to a cover where a singleton $x$ is duplicated into a pair $\langle
 x^1, x^2\rangle$ and $R(x^i, x^j)$ holds for $i,j\in \{1,2\}$ (let us still call $(X,R)$ such duplicating cover).  We let $Y$ to be the set of ordered 
 distinct pairs $\langle x_1,x_2\rangle$ from $X$ such that $R(x_1,x_2)$ holds in $(X,R)$. We let $R'(\langle x_1,x_2\rangle, \langle y_1,y_2\rangle)$ hold iff $\{x_1, x_2\}=\{y_1, y_2\}$. This turns $(Y, R')$ into a 1-step contact frame. The cover map $f: (Y, R') \longrightarrow (X,R)$ takes $\langle x_1,x_2\rangle$ to $x_1$.
\end{proof}

Now we are ready to show that Corollary~\ref{maincoro} applies.

\begin{theorem}\label{thm:con_amalg}
 The universal 
theory $\sf Con$ of contact algebras has the amalgamation property. Therefore, as it is also locally finite, $\mathsf{Con}$  has a model completion.
\end{theorem}

\begin{proof}
 As we observed in Section~\ref{sec:modelcompl}, it is sufficient to prove amalgamation for finite algebras (by local finiteness and by the
 compactness argument based on Robinson diagrams
 mentioned in Section~\ref{sec:modelcompl}).
 Finite algebras are dual to discrete Stone spaces, hence it is sufficient to show the following.
 
 \noindent
$(+)$~ Given finite nonempty sets $X_A, X_B, X_C$ endowed with reflexive and symmetric relations $R_A, R_B, R_C$ and given regular stable maps $f:(X_B, R_B) \to (X_A, R_A)$, $g:(X_B, R_B) \to (X_A, R_A)$, there exist $(X_D, R_D)$ (with reflexive and symmetric $R_D$) and regular stable maps $\pi_1:(X_D, R_D) \to (X_B, R_B)$, $\pi_2:(X_D, R_D) \to (X_C, R_C)$, such that $f\circ \pi_1=g\circ \pi_2$.
\par
 Statement $(+)$ is easily proved by taking as $(X_D, R_D), \pi_1, \pi_2$ the obvious pullback with the two projections.
\end{proof}

Theorem~\ref{thm:con_amalg} gives the possibility of applying Corollary~\ref{maincoro} to recognize admissible rules. We give here another algorithm, slightly different from that of Corollary~\ref{maincoro}.
We recall that ${\sf Con}^\star$ is the theory of existentially closed contact algebras~\cite{CK}.
The following result 
(given that ${\sf Con}$ is locally finite) 
is folklore (a detailed proof 
of the analogous statement for Brouwerian semilattices is 
in the ArXiv version of~\cite{CG19} as \cite[Proposition 2.16]{GhilardiCarai}).

\begin{theorem} \label{finite_extensions_contact}
Let $(\cB, \prec)$ be a contact algebra. We have that $(\cB, \prec)$ 
is existentially closed iff for any finite subalgebra $(\cB_0, \prec) \subseteq (\cB, \prec)$ and for any finite extension $(\cC, \prec) \supseteq (\cB_0, \prec)$ there exists an embedding $(\cC, \prec) \hookrightarrow{} (\cB, \prec)$
such that the following diagram commutes
\begin{center}
\begin{tikzcd}
(\cB_0, \prec) \arrow[r, hook] \arrow[d, hook] & (\cB, \prec) \\
(\cC, \prec) \arrow[ru, hook] & {}
\end{tikzcd}
\end{center}
\end{theorem}

\begin{example}
 Consider the $\Pi_2$-rule:
 $$
\leqno{(\rho9)} \quad \inference{(p \simpa p) \land (\varphi \simpa p) \land (p \simpa \psi)
  \to \chi}{(\varphi \simpa \psi) \to \chi}
$$
This rule is admissible in $\sf S^2IC$ \cite[Theorem 6.15]{BBSV}. We will now give an alternative and more automated proof of this result.
Translating $\Pi(\rho9)$ into the equivalent language of contact algebras, 
we obtain (see statement $(S9)$ from Section 6.3 of~\cite{BBSV})
\begin{equation}\label{eq:ex}
x \prec y \Rightarrow \exists z ~(z \prec z \fowedge x \prec z \prec y).
\end{equation}
According to Theorem~\ref{thm:mainadm}, we have to show that~\eqref{eq:ex} is provable in ${\sf Con}^\star$.
Note that~\eqref{eq:ex} expresses interesting (order-)topological  properties. It is valid on $(X, R)$ iff $R$ is a Priestley quasi-order \cite[Lemma 5.2]{BBSV16}.  Also it is valid on a compact Hausdorff space $X$ iff $X$ is a Stone space \cite[Lemma 4.11]{GB10}.

If we follow the procedure of Corollary~\ref{maincoro}
(which is based on Lemma~\ref{lem:qe}), we  first compute the quantifier-free formula equivalent in ${\sf Con}^\star$ to $\exists z ~(z \prec z \fowedge x \prec z \prec y)$ by taking the conjunction of the (finitely many) quantifier-free first-order formulas $\phi(x,y)$ which are implied (modulo ${\sf Con}$) by $z \prec z \fowedge x \prec z \prec y$: this is, up to equivalence, $x\prec y$. Now, in order to show the admissibility of $(\rho9)$ it is sufficient to observe that ${\sf Con}\models x \prec y\Rightarrow x \prec y$. 

As an alternative, we can rely on Theorem~\ref{finite_extensions_contact} and show that~\eqref{eq:ex} is true in every existentially closed contact algebra. To this aim, it is sufficient to enumerate all contact algebras $\cB_0$ generated by two elements $a,b$  such that $\cB_0\models a\prec b$ and to show that all such algebras embed in a contact algebra $\cC$ generated by three elements $a,b, c$ such that $\cC\models c \prec c \fowedge a \prec c \prec b$ (this can be done automatically for instance using a model finder tool). Both of the above procedures are heavy and not elegant, but they are nevertheless mechanical and do not require ingenious \emph{ad hoc} constructions (such as e.g., the construction of~\cite[Lemma 5.4]{BBSV}).
\end{example}

\subsection{Complexity Issues}\label{subsec:complexity} In this subsection, we will adopt the algorithm suggested by Corollary~\ref{maincoro} and Lemma~\ref{lem:qe} to get a co-\textsc{NExpTime} upper bound for deciding admissibility of $\Pi_2$-rules in 
$\sf S^2IC$.

In order to do that, we first need to study closer the satisfiability problem for quantifier-free formulas in the language of contact algebras. First notice that atomic formulas in such a language are all equivalent to formulas of the kind $t\prec u$, where $t,u$ are Boolean terms: this is because atoms of the kind $t=0$ are equivalent to $t\prec \neg t$, by the axioms of contact algebras. Second, we introduce a more manageable Kripke style equivalent semantics for satisfiability in finite contact algebras (since finitely generated contact algebras are finite and contact algebras axioms are all universal, to test satisfiability of a quantifier-free formula it is sufficient to inspect finite contact algebras).

A \emph{Kripke model} over a contact frame $(X,R)$ is a valuation $V:\mathsf{Prop} \longrightarrow \wp(X)$ from the set of propositional variables into the power set of $X$; we use the notation $\cM=(X,R,V)$ for such a Kripke model. For $x\in X$ and a Boolean formula (term) $F$, the notion $\cM\models_x F$ is defined inductively as follows: 
\begin{itemize}
 \item[-] $\cM\models_x p$ iff $x\in V(p)$, for atomic $p$;
 \item[-] $\cM\models_x F_1\wedge F_2$ iff ($\cM\models_x F_1$ and $\cM\models_x F_2$);
 \item[-] $\cM\models_x \neg F$ iff $\cM\not \models_x  F$.
\end{itemize}
For an atom $F\prec G$, we put $\cM\models F\prec G$ iff for all $x, y\in X$, we have that $\cM\models_x F$ and $R(x,y)$ imply
$\cM\models_y G$. Finally, for a quantifier-free formula $\phi$, the definition of $\cM\models \phi$ goes by induction as expected. The following lemma is clear.

\begin{lemma}
 Let $\cA$ be a finite contact algebra with dual contact frame $(X,R)$. For a quantifier-free formula $\phi$, we have that 
 $\phi$ is true in $\cA$ under some free variables assignment iff we have  $\cM\models \phi$ for some Kripke model
 $\cM=(X,R,V)$.
\end{lemma}

In our context, covers play the same role as p-morphisms in modal logic. A \emph{cover} of a Kripke model  $\cM=(X,R,V)$ is a Kripke model  $\cM'=(X',R',V')$  together with a regular stable map 
 $f: (X',R')\longrightarrow (X,R)$  of the underlying contact frames such that for every propositional variable $p$, we have 
that $V'(p)= f^{-1}(V(p))$.

\begin{lemma}\label{lem:p}
  Let $\cM'=(X',R',V')$ be a cover of $\cM=(X,R,V)$ $($via a suitable $f)$. Then for every quantifier-free formula $\phi$, we have that $\cM\models \phi$ iff $\cM'\models \phi$.
\end{lemma}

\begin{proof}
 This follows from the fact that $f^{-1}$ induces, as we know, an embedding of the contact algebras dual to $(X,R)$ and $(X',R')$. 
\end{proof}

\begin{lemma}\label{lem:NP}
 A quantifier-free formula $\phi$ in the language of contact algebras is satisfiable iff it is satisfiable in  a finite quadratic size 1-step contact frame. Thus the satisfiability problem for $\phi$ is in NP.
\end{lemma}

\begin{proof}
 First observe that 
$\phi$ is satisfiable iff there is a consistent Boolean assignment to the atoms of $\phi$ satisfying $\phi$ from the point of view of propositional logic. To show that a  candidate Boolean assignment is satisfiable one translates positive atoms (as well as the reflexivity and symmetry  conditions for the relation of a contact frame) into universally quantified Horn clauses in first-order logic using at most two universally quantified variables. 
Negative atoms $F\not\prec G$ translate into $\exists x \exists y \;(F(x)\fowedge R(x,y) \fowedge \foneg G(y))$; skolemization of these literals introduces two Skolem constants for each of them.  Thus, the overall universal Horn formula to be checked for satisfiability  has a finite Herbrand universe of linear size. 

Since the Herbrand universe is of linear size,
 $\phi$ is satisfiable iff it is satisfiable in a  linear size finite contact frame (alternative ways to prove this  arise from translations into $\mathsf{S5}_U$, see~\cite{BBSV}). Then the fact that 1-step quadratic contact frames suffice follows from Lemma~\ref{lem:1step}.
\end{proof}

According to Theorem~\ref{thm:mainadm}, the rule $\rho$ is \emph{not} admissible iff
$T_{S^2I}^\star\not\models \Pi(\rho)$, where 
$T_{S^2I}^\star$
is the model completion of the theory of 
simple symmetric strict implication algebras. Since we want to go through the 
equivalent theory given by the model completion 
${\sf Con}^\star$
of the theory of contact algebras, 
in view of the  Lemma~\ref{lem:NP}, to get our co-\textsc{NExpTime} upper bound, it is sufficient to prove that \emph{the computation of the quantifier-free formula 
$\phi^\star(\ux)$
equivalent in ${\sf Con}^\star$ to $\exists \uy \phi(\ux, \uy)$}  (for any quantifier-free formula $\phi(\ux, \uy)$) \emph{is exponentially large and can be computed in exponential time} (because then Lemma~\ref{lem:NP} would apply). 
However, 
since there are  double exponentially many non-equivalent quantifier-free formulas built up from a finite set of variables in the language of contact algebras, this  is not obvious. The situation is similar to the problem of showing an exponential bound for  the computation of uniform interpolants in $\mathsf{S5}$~\cite{GLWZ06} and in fact we will solve it 
by adapting the technique of~\cite{GLWZ06} to our context.

Let $N$ be the number of distinct atoms (i.e.,  atomic formulas) occurring in  $\phi(\ux, \uy)$ and let us consider the Kripke models built up on finite 1-step contact frames having at most 
$2N$
elements (they are exponentially many). Partition them into classes $K_1,\dots, K_m$ in such a way that two  models are in the same class $K_i$ iff they satisfy the same atoms from $\phi$. To every $\cM\in K_i$ associate the formula
\begin{equation}\label{eq:CM}
 \chi(\cM) := \bigwedge \{ t(\ux)\not\prec u(\ux) ~\mid~ t(\ux), u(\ux) ~{\rm are~Boolean~terms~}~{\rm s.t.}~ \cM \not
 \models t(\ux)\prec u(\ux)\}.
\end{equation}
Let also 
\begin{equation}\label{eq:CKi}
\theta_i := \bigwedge \{ t(\ux)\prec u(\ux) ~\mid~ t(\ux), u(\ux) ~{\rm are~Boolean~terms~}~~{\rm s.t.}~ \cM' \models
t(\ux)\prec u(\ux)~{\rm for~all~}\cM'\in K_i\}.
\end{equation}
We claim that the formula we need is 
\begin{equation}\label{eq:C}
\phi^\star(\ux)~:=~\bigvee_{i=1}^m ( \theta_i \fowedge \bigvee_{\cM\in Ki} \chi(\cM) ).
\end{equation}
Notice that this is (simply) exponential.

According to Lemma~\ref{lem:qe}, we must show that $\phi\Rightarrow \phi^\star$ holds in $T_{\cS}$ and that if $\psi(\ux)$ is such that $T_{\cS}\models \phi\Rightarrow \psi$, then $T_{\cS}\models \phi^\star\Rightarrow \psi$; by Lemma~\ref{lem:NP}, all validity tests can be performed in Kripke models over finite contact 1-step frames.

First consider a Kripke model $\cN=(X,R,V)$ based on a 1-step contact frame such that $\cN\models \phi$. Restrict the model to a submodel by picking one witness pair $x,y$ for every atom $u_1(\ux,\uy)\prec u_2(\ux,\uy)$ such that $\cN\models_x u_1$, $\cN\not\models_y u_2$ and $R(x,y)$.  The Kripke model $\cM$ obtained by this restriction is such that $\cM\models \phi$, it has the size at most $2N$ and it thus belongs (up to isomorphism) to a certain partition $K_i$. Clearly we have $\cN\models \chi(\cM)$. We also have $\cN\models \theta_i$ because for every $x',y'\in X$ such that $R(x',y')$, we can always pick witness points so as to build a submodel $\cM'$ of $\cN$ belonging (up to isomorphism) to $K_i$ and including $(x',y')$. Thus we obtain $\cN\models \phi^\star$, as desired.

Suppose now that $T_{\cS}\not\models \phi^\star\Rightarrow \psi$, i.e., there is a  Kripke model based on a 1-step contact frame such that $\cN\models \phi^\star\fowedge \foneg\psi$. In order to show that $T_{\cS}\not\models \phi\Rightarrow \psi$ we proceed as follows.
Let $\cN$ be $(X,R,V)$. 
We build $\cN'=(X',R',V')$ and a regular stable map $f:(X',R')\longrightarrow (X,R)$ in such a way that $\cN'\models \phi(\ux, \uy)$ and $V'(p)=f^{-1}(V(p))$ for all $p\in \ux$ (this guarantees that
$\cN'\not \models \psi(\ux)$, 
by Lemma~\ref{lem:p}).

Since $\cN\models \phi^\star(\ux)$, there are $i$ and $\cM\in K_i$ such that
$\cN \models \theta_i \fowedge  \chi(\cM)$.
Recall that $\cM$ is based on a 1-step finite contact frame and suppose that $\cM=(X_0, R_0, V_0)$.
For every $x\in X_0$ we can build the atoms 
$$
t^{+}_x=\bigwedge \{p \in \ux \mid x\in V_0(p)\}\fowedge \foneg\bigvee \{p \in \ux \mid x\not\in V_0(p)\} 
$$ 
and
$$
t^{-}_x=\foneg\, t^{+}_x.
$$
Now notice that for every distinct pair $\langle x_1,x_2\rangle$
such that $R_0(x_1,x_2)$ holds in $\cM$, the atom $t_{x_1}^{+}\prec t_{x_2}^{-}$ is false in $\cM$ precisely because of the pair $(x_1, x_2)$. Since $\cN \models  \chi(\cM)$ there must be a pair (not necessarily formed by distinct 
elements) $\langle w_1, w_2\rangle$ in $\cN$ such that $R(w_1,w_2)$ and 
$\cN\models_{w_1} t^{+}_{x_1}$, $\cN\not\models_{w_2} t^{-}_{x_2}$, 
which means that for evey $p\in \ux$, we have 
$x_1\in V_0(p) \Leftrightarrow w_1\in V(p)$ and  $x_2\in V_0(p) \Leftrightarrow w_2\in V(p)$. Since $\cM$ is a 1-step contact frame, this defines a stable map 
$f_0: (X_0, R_0)\longrightarrow (X,R)$ which preserves the satisfiability of the variables $\ux$. However, this map may  not be regular (not even surjective), to make it regular we need a further simple adjustement: we take as $(X',R')$ the disjoint union of $(X, R)$ with $(X_0,R_0)$ and as $f$ the identity map coupled with $f_0$. 
This obviously gives a regular stable map.
It remains to define the forcing $V'$ on $(X',R')$ for the variables $\uy$. This must be done  in such a way that $\phi(\ux,\uy)$ becomes true. 

For $q\in \uy$ and 
$x$ in the $X_0$-part of $X'$ we just use the satisfiability in $\cM$, that is we let $x\in V'(q)$ hold iff $x\in V_0(q)$.
Let now consider a distinct pair $\langle x, y\rangle$ in the $X$-part of $X'$ such that $R(x,y)$ holds. We have that 
$\cN\not\models t_{x}^{+}\prec t_{y}^{-}$ and, since $\cN\models \theta_i$, there must be some $\cM_{xy}\in K_i$ with 
$\cM_{xy} \not\models t_{x}^{+}\prec t_{y}^{-}$. This means that for some $x',y'$ in the support of $\cM_{xy}$ we have 
$\cM_{xy}\models_{x'} t_{x}^{+}$ and $\cM_{xy}\models_{y'} t_{y}^{-}$ (which is the same as $\cM_{xy}\models_{y'} t_{y}^{+}$). For $q\in \uy$, we let $V'(q)$ contain $x$ (resp.~$y$) iff we have $\cM_{xy}\models_{x'} q$ (resp.~$\cM_{xy}\models_{y'} q$). Notice that the same relations holds for $q\in \ux$ because $\cM_{xy}\models_{x'} t_{x}^{+}$ and 
$\cM_{xy}\models_{y'} t_{y}^{-}$. Now $\phi(\ux,\uy)$ holds in $\cN'=(X',R',V')$ because exactly the same atoms from $\phi$ satisfied in all members  of the class $K_i$ are true in 
$\cN'$.

We have therefore proved the following result.

\begin{theorem}\label{thm:complexity}
 The problem of recognizing the admissibility of a $\Pi_2$-rule in the symmetric strict implication calculus $\mathsf{S^2IC}$  
 is co-\textsc{NExpTime}-complete.
\end{theorem}

\begin{proof}
According to Theorem~\ref{thm:mainadm}, the $\Pi_2$-rule $\rho$ given in Definition~\ref{defirules} is \emph{not} admissible in $\mathsf{S^2IC}$ iff in the model completion $T_{S^2I}^\star$ of the theory of simple symmetric strict implication algebras, the formula
\[
\Pi(\rho) := \forall \ux\,\forall z\, \exists \uy  \Big( G(\ux) \nleq z 
\Rightarrow
 F(\ux,\uy) \nleq z \Big)
\]
is not provable. To check this, we should eliminate the existential quantifiers 
from $\exists \uy  ( G(\ux) \nleq z 
\Rightarrow
 F(\ux,\uy) \nleq z) $ in $T_{S^2I}^\star$, then get a universal formula 
 $\forall \ux\,\forall z\, \psi(\ux,z)$, and finally check $\neg\psi(\ux,z)$ for satisfiability in $T_{S^2I}^\star$ (or, which is the same, in $T_{S^2I}$). 
 In view of Proposition~\ref{prop:translation} and Lemma~\ref{lem:qe}, we can equivalently apply these operations in 
 ${\sf Con}^\star/\sf Con$
 to the translation $\tau_2$ of $G(\ux) \nleq z 
\Rightarrow
 F(\ux,\uy) \nleq z$. 
 
 In principle, $\tau_2$  may cause an exponential blow-up in the third step of its computation, but since our first task is  to  eliminate the existential quantifiers from  $\exists \uy  \;\tau_2( G(\ux) \nleq z 
\Rightarrow
 F(\ux,\uy) \nleq z) $, we can just eliminate the existential quantifiers from the equivalent formula~\eqref{eq:firsttranslation} obtained in the second step of the computation of $\tau_2$: such a formula is only linearly long, and consequently, as explained above, our quantifier elimination procedure takes exponential time and produces an exponentially long formula. Thus, in the end, Lemma~\ref{lem:NP} gives our desired \textsc{NExpTime} upper bound.

For the lower bound, we notice that 
 in~\cite{GLWZ06} it is shown that checking conservativity in $\mathsf{S5}$ is co-\textsc{NExpTime}-complete. Conservativity is trivially translated into 
 admissibility of  $\Pi_2$-rules for logics like $\mathsf{S5}$ enjoying interpolation (see Theorem~\ref{thm:cons}) and on the other hand $\mathsf{S5}$ is a subsystem of $\mathsf{S^2IC}$. Thus, it is sufficient to show that a $\Pi_2$-rule in the restricted language of $\mathsf{S5}$ which is admissible in $\mathsf{S5}$ is also admissible in $\mathsf{S^2IC}$ (the vice versa is obvious). To this aim, we apply the admissibility criterion given by Theorem~\ref{thm:criterion}. 
 Consider a  $\Pi_2$-rule $\rho$ as given in Definition~\ref{defirules}, which is in the language of $\mathsf{S5}$ and is admissible in $\mathsf{S5}$.
 Let  $\cB$ be a simple $S^2I$-algebra; according to Theorem~\ref{thm:criterion}, its $\mathsf{S5}$-reduct (which is nothing but a Boolean algebra, being the algebra simple) embeds into an algebra $\cB'$ satisfying $\Pi(\rho)$. Thus, it is sufficient to apply the Lemma below (we exploit once again the equivalence between contact algebras and simple symmetric strict implication algebras). 
 \end{proof}
 
 \begin{lemma}
  Given a contact algebra $\cB$ and a Boolean algebra $\cB'$ extending it, it is possible to give $\cB'$ a structure of a contact algebra in such a way that the embedding preserves also the contact algebra structure.
 \end{lemma}

\begin{proof}
 We prove the dual statement using Theorem~\ref{thm:duality}. Let $Y$ be a Stone space, $(X,R)$ be an object of $\mathsf{StR}$ and let $f: Y\longrightarrow X$ be a continuous surjective map. We endow $Y$ with the relation $\tilde R$ given by $y_1 \tilde R y_2$ iff $f(x_1) R f(y_2)$. 
 Since $\tilde R$ is closed, it turns out that $f:(Y,\tilde R) \longrightarrow (X,R)$ is a morphism in $\mathsf{StR}_e$.
\end{proof}

\section{Finite axiomatization of ${\sf Con}^\star$}\label{subsec:axiomatization} 

Theorem~\ref{finite_extensions_contact} implicitly supplies an infinite set of 
 axioms for ${\sf Con}^\star$, the model completion of the theory of contact algebras. This axiomatization is not however very informative, as it follows from generic model-theoretic facts. In this section, we supply a better axiomatization following
 a strategy similar to the one used in~\cite{DJ18} for the case of amalgamable locally finite varieties of Heyting algebras and  in~\cite{CG19} for the case of  Brouwerian semilattices.
This axiomatization is finite and is described by the following theorem, which is the main result of this section.

\begin{theorem}\label{thm:axioms_model_completion}
An axiomatization of $\Con^\star$ is given by the axioms of contact algebra together with the following sentences:
\begin{align*}
&\forall a,b_1,b_2  \ (a \neq 0 \fowedge (b_1 \vee b_2)\wedge a=0 \fowedge a \prec a \vee b_1 \vee b_2 \Rightarrow \exists a_1,a_2 \ (a_1 \vee a_2=a \fowedge a_1 \wedge a_2=0 \tag{s1} \\
& \hspace{1.7cm} \fowedge a_1 \neq 0 \fowedge a_2 \neq 0 \fowedge a_1 \prec a_1 \vee b_1 \fowedge a_2 \prec a_2 \vee b_2  )), \\[0.2cm]
& \forall a,b \ (a \wedge b =0 \fowedge a \not\prec \neg b \Rightarrow \exists a_1, a_2 \ (a_1 \vee a_2=a \fowedge a_1 \wedge a_2 = 0 \fowedge a_1 \not\prec \neg b \fowedge a_2 \not\prec \neg b \fowedge a_1 \prec \neg a_2)),  \tag{s2}\\[0.2cm]
& \forall a \ (a \neq 0 \Rightarrow \exists a_1,a_2 \ (a_1 \vee a_2 = a \fowedge a_1 \wedge a_2=0 \fowedge a_1 \prec a \fowedge a_1 \not\prec a_1)). \tag{s3}
\end{align*}
\end{theorem}

Notice that the axioms (s1), (s2), (s3) are similar to the splitting axioms of the axiomatizations appearing in \cite{DJ18} and \cite{CG19}.
We will prove Theorem~\ref{thm:axioms_model_completion} by employing Theorem~\ref{finite_extensions_contact} and the duality between $\mathsf{Con}_e$ and $\mathsf{StR}_e$ to characterize the duals of existentially closed algebras. We first show that it is enough to work with finite minimal extensions.

\begin{definition}
If $(\cC, \prec)$ is a contact algebra extending the contact algebra $(\cB_0, \prec)$, we say that such an extension is \emph{minimal} if it is proper and it does not contain any other proper extension of $(\cB_0, \prec)$.
\end{definition}

Using the Duality Theorem~\ref{thm:duality} restricted to the finite discrete case, we can characterize the dual spaces $(X_{\cC},R_{\cC})$ and $(X_{\cB_0}, R_{\cB_0})$ and the dual stable map $f:(X_{\cC},R_{\cC})\to (X_{\cB_0}, R_{\cB_0})$ corresponding to finite minimal extensions.
 \begin{proposition}\label{prop:data} Let $(\cB_0,\prec)\hookrightarrow (\cC, \prec)$ be an embedding between finite contact algebras, with dual regular stable map $f:(X_{\cC},R_{\cC})\to (X_{\cB_0}, R_{\cB_0})$. The embedding is minimal iff  $($up to isomorphism$)$
 there are a finite set $Y$, finite subsets $S_1, S_2\subseteq Y$ and elements $x\in X_{\cB_0}, x_1\in X_{\cC}, x_2\in X_{\cC}$ such that: 
 \begin{description} 
  \item[{\rm (i)}] $X_{\cB_0}$ is the disjoint union $Y\oplus\{x\}$;
  \item[{\rm (ii)}] $X_{\cC}$ is the disjoint union $Y\oplus\{x_1, x_2\}$;
  \item[{\rm (iii)}] $f$ restricted to $Y$ is the identity map and $f(x_1)=f(x_2)=x$;
  \item[{\rm (iv)}] the restrictions of $R_{\cC}$ and of $R_{\cB_0}$ to $Y$ coincide;
  \item[{\rm (v)}] $R_{\cC}[x_1]\setminus\{x_1\}=S_1$ and $R_{\cC}[x_2]\setminus\{x_2\}=S_2$;
  \item[{\rm (vi)}] $R_{\cB_0}[x]\setminus\{x\}=S_1\cup S_2$.
 \end{description}
\end{proposition}

\begin{proof}
 First notice that, as a consequence of~\eqref{eq:emb},  if the cardinality of $X_{\cB_0}$ and of $X_{\cC}$ is the same, then $f$ is an isomorphism. This is seen as follows: we already observed  that condition~\eqref{eq:emb} implies surjectivity and in case of the same finite cardinality surjectivity implies injectivity. Preservation and reflection of the relation follow by stability and~\eqref{eq:emb} again. 
 
 In addition, if the cardinality of $X_{\cC}$ is equal to the cardinality of $X_{\cB_0}$ plus one (this is precisely the case mentioned in the statement of the proposition), then $f$ cannot be properly factored, hence it is minimal. We show that all minimal maps arise in this way.
 
 In general, if the cardinality of $X_{\cC}$ is bigger than the cardinality of $X_{\cB_0}$, we can define the following factorization of $f$. Pick some $x\in X_{\cB_0}$ having more than one preimage and split 
 $f^{-1}(\{x\})$ as $T_1\cup T_2$, where $T_1, T_2$ are disjoint and non-empty. We have that $X_{\cC}$ is the disjoint union $X\oplus T_1\oplus T_2$ for some set $X$ and $X_{\cB_0}$ is the disjoint union $Y\oplus\{x\}$ for some set $Y$. Define a discrete dual space $(Z, R_Z)$ as follows. $Z$ is the disjoint union $Y\oplus\{x_1, x_2\}$ for new $x_1, x_2$ and $R_Z$ is the reflexive and symmetric closure of the following sets of pairs: (i) the pairs $(z_1, z_2)$ for $z_1 R_{\cB_0} z_2$ and $z_1, z_2\in Y$; (ii) the pairs $(x_i, u)$ for 
 $u\in f(R_{\cC}[T_i])$ ($i=1,2$); (iii) the pair $(x_1,x_2)$, but only  in case $T_1\cap R_{\cC}[T_2]\not=\emptyset$. Then it is easily seen that $f$ factorizes as $h\circ \tilde f$ in $\mathsf{StR}_e$, where: (I) $\tilde f$ maps 
 $T_1$ to $x_1$, $T_2$ to $x_2$ and acts as $f$ on $X$; (II) $h$ is the identity on $Y$ and maps both $x_1, x_2$ to $x$.
 
 Now $h$ produces the  data required by the proposition  and $\tilde f$ must be an isomorphism if $f$ is minimal.
\end{proof}

\begin{remark}
\mbox{}\begin{enumerate}
\item The conditions (i)--(vi) in Proposition~\ref{prop:data} determine uniquely the finite minimal extension over the contact algebras dual to $(X_{\cB_0}, R_{\cB_0})$  except for a detail: they do not specify whether we have $x_1R_{\cC} x_2$ or not. So the data $x,S_1, S_2$ and  $Y=X_{\cB_0}\setminus\{x\}$ (lying \emph{inside $X_{\cB_0}$}) determine in fact \emph{two} minimal extensions of the contact algebra dual to $(X_{\cB_0}, R_{\cB_0})$.
\item It is an immediate consequence of the proof of Proposition~\ref{prop:data} that every finite extension of contact algebras can be decomposed into a finite chain of finite minimal extensions. Thus, Theorem~\ref{finite_extensions_contact} still holds if we limit its statement to finite \emph{minimal} extensions.
\end{enumerate}
\end{remark}

Thus, by dualizing Theorem~\ref{finite_extensions_contact}, we obtain the following characterization of the contact frames that are dual to existentially closed contact algebras.

\begin{proposition}\label{prop:dual_exclosed_via_extensions}
The contact frame $(X,R)$ is dual to an existentially closed contact algebra iff for every finite contact frame $(Y_0,R_0)$, every regular stable map $f: (Y_1,R_1) \to (Y_0,R_0)$ dual to a finite minimal extension of contact algebras, and every regular continuous stable map $g:(X,R) \to (Y_0,R_0)$ there exists a continuous regular stable map $h:(X,R) \to (Y_1,R_1)$ such that $f \circ h =g$.
\[
\begin{tikzcd}
(Y_0,R_0) & (X,R) \arrow[l, two heads, "g"'] \arrow[dl, two heads, "h"] \\
(Y_1,R_1)  \arrow[u, two heads, "f"]  &
\end{tikzcd}
\]
\end{proposition}

We reformulate this characterization of duals of existentially closed contact algebras in terms of partitions.

\begin{lemma}\label{lem:ex closed and partitions}
Let $(X,R) \in \Str$. The contact algebra $(X,R)^\ast$ is existentially closed iff for each finite partition $\Pcal$ of $X$ into clopens, $A \in \Pcal$, and $\S_1, \S_2 \subseteq \Pcal$ with $\S_1 \cup \S_2=\{ C \in \Pcal \setminus \{A\} \mid A \cap R[C] \neq \emptyset \}$, there exist two nonempty clopens $A_1, A_2$  such that $A_1 \cup A_2 =A$, $A_1 \cap A_2=\emptyset$ and for each $C \in \Pcal \setminus \{A\}$
\[
A_1 \cap R[C] \neq \emptyset  \ \ \mbox{iff} \ \ C \in \S_1, \qquad A_2 \cap R[C] \neq \emptyset \ \ \mbox{iff} \ \ C \in \S_2 \qquad \mbox{and} \qquad A_1 \cap R[A_2] = \emptyset.
\]
and there exist two nonempty clopens $A_1', A_2'$ such that $A_1' \cup A_2' =A$, $A_1' \cap A_2'=\emptyset$ and for each $C \in \Pcal \setminus \{A\}$
\[
A_1' \cap R[C] \neq \emptyset \ \ \mbox{iff} \ \ C \in \S_1, \qquad A_2' \cap R[C] \neq \emptyset \ \ \mbox{iff} \ \ C \in \S_2 \qquad \mbox{and} \qquad A_1' \cap R[A_2'] \neq \emptyset.
\]
\end{lemma}

\begin{proof}
This is a consequence of Propositions~\ref{prop:data} and~\ref{prop:dual_exclosed_via_extensions}, and the fact that continuous regular stable maps from $(X,R) \in \Str$ into finite objects of $\Str$ correspond to finite partitions of $X$ into clopens. Indeed, a continuous regular stable map $f:(X,R) \to (Y,R')$ onto a finite contact frame induces the partition $\Pcal=\{ f^{-1}(y) \mid y \in Y \}$. On the other hand, if $\Pcal$ is a finite partition of $X$ into clopens, the quotient map $f:(X,R) \to (\Pcal,R_\Pcal)$ is a continuous regular stable map, where $A \relr_\Pcal B$ iff $A \cap R[B] \neq \emptyset$ for any $A,B \in \Pcal$.
\end{proof}

We are ready to show that the following conditions, which dually correspond to the axioms (s1), (s2), (s3) of Theorem~\ref{thm:axioms_model_completion}, characterize the contact frames $(X,R)$ dual to existentially closed contact algebras.
\begin{enumerate}
\item[(S1)] If $A,B_1,B_2$ are clopens of $X$ with 
\[
A \neq \emptyset, \quad (B_1 \cup B_2) \cap A = \emptyset \quad \mbox{and} \quad R[A] \subseteq A \cup B_1 \cup B_2,
\]
then there exist $A_1,A_2$ clopens of $X$ such that 
\begin{align*}
A_1 \cup A_2 =A, \quad A_1 \cap A_2=\emptyset, \quad A_1 \neq \emptyset, \quad A_2 \neq \emptyset, \quad R[A_1] \subseteq A_1 \cup B_1, \quad \mbox{and} \quad R[A_2] \subseteq A_2 \cup B_2.
\end{align*}
\item[(S2)] If $A,B$ are clopens of $X$ with
\[
A \cap B= \emptyset, \quad A \cap R[B] \neq \emptyset,
\]
then there exist $A_1,A_2$ clopens of $X$ such that 
\begin{equation*}
A_1 \cup A_2 =A, \quad A_1 \cap A_2=\emptyset, \quad
A_1 \cap R[B] \neq \emptyset, \quad A_2 \cap R[B] \neq \emptyset, \quad \mbox{and} \quad A_1 \cap R[A_2] = \emptyset.
\end{equation*}
\item[(S3)] If $A$ is a nonempty clopen of $X$, then there exist $A_1,A_2$ clopens of $X$ such that
\begin{align*}
A_1 \cup A_2 =A, \quad A_1 \cap A_2=\emptyset, \quad R[A_1] \subseteq A, \quad \mbox{and} \quad R[A_1] \nsubseteq A_1.
\end{align*}
\end{enumerate}

\begin{lemma}\label{lem:exclosedimpliesS1S2S3}
Let $(X,R) \in \Str$. If the contact algebra $(X,R)^\ast$ is existentially closed, then \emph{(S1)}, \emph{(S2)}, and \emph{(S3)} hold in $(X,R)$.
\end{lemma}

\begin{proof}
(S1) Let $A,B_1,B_2$ be clopens of $X$ such that $A \neq \emptyset$, $(B_1 \cup B_2) \cap A = \emptyset$, and $R[A] \subseteq A \cup B_1 \cup B_2$. Let $\Pcal$ be the partition obtained from
\[
\{A,B_1 \setminus B_2, B_2 \setminus B_1, B_1 \cap B_2, X \setminus (A \cup B_1 \cup B_2) \}
\]
after possibly removing the empty set from its elements.
Let $\S_i=\{ C \in \Pcal  \setminus \{A\}  \mid A \cap R[C] \neq \emptyset, \  C \subseteq B_i\}$ for $i=1,2$. Since $R[A] \subseteq A \cup B_1 \cup B_2$, we have $\S_1 \cup \S_2 = \{ C \in \Pcal \setminus \{A\} \mid A \cap R[C] \neq \emptyset \}$. Therefore, by Lemma~\ref{lem:ex closed and partitions} there exist $A_1, A_2$ nonempty clopens such that $A_1 \cup A_2 =A$, $A_1 \cap A_2=\emptyset$ and for each $C \in \Pcal \setminus \{A\}$:
\[
A_1 \cap R[C] \neq \emptyset \ \mbox{iff} \ C \in \S_1, \qquad A_2 \cap R[C] \neq \emptyset \ \mbox{iff} \ C \in \S_2 \qquad \mbox{and} \qquad A_1 \cap R[A_2] = \emptyset.
\]
It follows that $R[A_1] \subseteq A_1 \cup B_1$ and $R[A_2] \subseteq A_2 \cup B_2$.

(S2) Let $A,B$ be clopens of $X$ such that $A \cap B= \emptyset$ and $A \cap R[B] \neq \emptyset$. Since $A \cap R[B] \neq \emptyset$, both $A$ and $B$ are not empty. Let $\Pcal$ be the partition obtained from
\[
\{A, B, X \setminus (A \cup B) \}
\]
after possibly removing $X \setminus (A \cup B) $ if it is empty. Let $\S_1=\S_2=\{ C \in \Pcal \setminus \{A\} \mid A \cap R[C] \neq \emptyset \}$. By Lemma~\ref{lem:ex closed and partitions} there exist $A_1, A_2$ nonempty clopens such that $A_1 \cup A_2 =A$, $A_1 \cap A_2=\emptyset$ and for each $C \in \Pcal \setminus \{A\}$:
\[
A_1 \cap R[C] \neq \emptyset \ \mbox{iff} \ A \cap R[C] \neq \emptyset, \qquad A_2 \cap R[C] \neq \emptyset \ \mbox{iff} \ A \cap R[C] \neq \emptyset \qquad \mbox{and} \qquad A_1 \cap R[A_2] = \emptyset.
\]
It follows that $A_1 \cap R[B] \neq \emptyset$ and $A_2 \cap R[B] \neq \emptyset$.

(S3) Let $A$ be a nonempty clopen of $X$. Let $\Pcal$ be the partition obtained from
\[
\{A, X \setminus A \}
\]
after possibly removing $X \setminus A$ if it is empty. Let $\S_1=\emptyset$ and $\S_2=\{ C \in \Pcal \setminus \{A\} \mid A \cap R[C] \neq \emptyset \}$. Therefore, by Lemma~\ref{lem:ex closed and partitions} there exist $A_1, A_2$ nonempty clopens such that $A_1 \cup A_2 =A$, $A_1 \cap A_2=\emptyset$ and for each $C \in \Pcal \setminus \{A\}$:
\[
A_1 \cap R[C] = \emptyset, \qquad A_2 \cap R[C] \neq \emptyset \ \mbox{iff} \ A \cap R[C] \neq \emptyset \qquad \mbox{and} \qquad A_1 \cap R[A_2] \neq \emptyset.
\]
Therefore, $R[A_1] \subseteq A$ and $R[A_1]  \nsubseteq A_1$.
\end{proof}

\begin{lemma}\label{lem:S1S2S3implyexclosed}
Let $(X,R) \in \Str$. If \emph{(S1)}, \emph{(S2)}, and \emph{(S3)} hold in $(X,R)$, then the contact algebra $(X,R)^\ast$ is existentially closed.
\end{lemma}

\begin{proof}
We will show using Lemma~\ref{lem:ex closed and partitions} that if (S1), (S2), and (S3) hold in $(X,R)$, then $(X,R)^\ast$ is existentially closed. Let $\Pcal$ be a finite partition of $X$ into nonempty clopens, $A \in \Pcal$, $\S_1, \S_2 \subseteq \Pcal$ such that $\S_1 \cup \S_2=\{ C \in \Pcal \setminus \{A\} \mid A \cap R[C] \neq \emptyset \}$. Let $\S_1 \cup \S_2 = \{ B_1, \ldots, B_n \}$. 

First we consider the case when $\S_1$ or $\S_2$ is empty. We can assume without loss of generality that $\S_1=\emptyset$ and hence that $\S_2=\{ C \in \Pcal \setminus \{A\} \mid A \cap R[C] \neq \emptyset \}$. Apply (S1) to $A, \emptyset, B_1 \cup \cdots \cup B_n$ to get $A_1,A_2$ clopens such that 
\begin{align*}
& A_1 \cup A_2 =A, \quad A_1 \cap A_2=\emptyset, \quad A_1 \neq \emptyset, \quad A_2 \neq \emptyset,\\
& R[A_1] \subseteq A_1, \quad R[A_2] \subseteq A_2 \cup B_1 \cup \cdots \cup B_n.
\end{align*}
Therefore, for each $C \in \Pcal \setminus \{A\}$
\[
A_1 \cap R[C] = \emptyset, \qquad A_2 \cap R[C] \neq \emptyset \ \ \mbox{iff} \ \ C \in \S_2 \qquad \mbox{and} \qquad A_1 \cap R[A_2] = \emptyset.
\]

Now assume that both $\S_1$ and $\S_2$ are not empty. We want to split $A$ into $n$ disjoint nonempty clopens $E_1, \ldots, E_n$. 
If $n=1$, let $E_1=A$. If $n >1$, apply (S1) to $A, B_1, B_2 \cup \cdots \cup B_n$ to get $D_{1, 1},D_{1, 2}$ clopens such that 
\begin{align*}
& D_{1, 1} \cup D_{1, 2} =A, \quad D_{1, 1} \cap D_{1, 2}=\emptyset, \quad D_{1, 1} \neq \emptyset, \quad D_{1, 2} \neq \emptyset,\\
& R[D_{1, 1}] \subseteq D_{1, 1}, \quad R[D_{1, 2}] \subseteq D_{1, 2} \cup B_2 \cup \cdots \cup B_n.
\end{align*}
Then we define recursively $D_{i, 1}, D_{i, 2}$ for each $i=2, \ldots, n-1$  by applying (S1) to $D_{i-1, 2}, B_i, B_{i+1} \cup \cdots \cup B_n$. Thus, we have that
\begin{align*}
& D_{i, 1} \cup D_{i, 2} =D_{i-1, 2}, \quad D_{i, 1} \cap D_{i, 2}=\emptyset, \quad D_{i, 1} \neq \emptyset, \quad D_{i, 2} \neq \emptyset,\\
& R[D_{i, 1}] \subseteq D_{i, 1} \cup B_i, \quad R[D_{i, 2}] \subseteq D_{i, 2} \cup B_{i+1} \cup \cdots \cup B_n.
\end{align*}
Let $E_i=D_{i,1}$ for $i=1, \ldots, n-1$ and $E_n=D_{n-1, 2}$.
This yields a family of nonempty clopens $E_1, \ldots, E_n$ such that 
\begin{align*}
E_1 \cup \cdots \cup E_n=A \ \mbox{ and } \ E_i \cap R[E_j]= \emptyset \mbox{ if } i \neq j
\end{align*}
and for each $C \in \Pcal \setminus \{ A \}$:
\[
E_i \cap R[C] \neq \emptyset \ \mbox{iff} \ C=B_i.
\]
The next step consists of splitting $E_i$ into two disjoint clopens for each $i$ such that $B_i \in \S_1 \cap \S_2$. Apply (S2) to $E_i,B_i$ for each $i=1, \ldots, n$ such that $B_i \in \S_1 \cap \S_2$. Thus, there exist $E_{i,1}, E_{i,2}$ clopens of $X$ such that 
\begin{align*}
E_{i,1} \cup E_{i,2} =E_i, \quad E_{i,1} \cap E_{i,2} =\emptyset,  \quad E_{i,1} \cap R[B_i] \neq \emptyset, \quad E_{i,2} \cap R[B_i] \neq \emptyset, \quad E_{i,1} \cap R[E_{i,2}] = \emptyset.
\end{align*}
We are finally ready to define $A_1$ and $A_2$. Let
\begin{align*}
A_1 &=\bigcup \{ E_i \mid B_i \in \S_1 \setminus \S_2 \} \cup \bigcup \{E_{i,1} \mid B_i \in \S_1 \cap \S_2 \},\\
A_2 &=\bigcup \{ E_i \mid B_i \in \S_2 \setminus \S_1 \} \cup \bigcup \{E_{i,2} \mid B_i \in \S_1 \cap \S_2 \}.
\end{align*}
It follows that $A_1,A_2$ are nonempty clopens such that $A_1 \cup A_2 =A$, $A_1 \cap A_2=\emptyset$, and for each $C \in \Pcal \setminus \{A\}$:
\[
A_1 \cap R[C] \neq \emptyset \ \mbox{iff} \ C \in \S_1, \qquad A_2 \cap R[C] \neq \emptyset \ \mbox{iff} \ C \in \S_2 \qquad \mbox{and} \qquad A_1 \cap R[A_2] = \emptyset.
\]
It remains to construct $A_1',A_2'$. Apply (S3) to $A_2$ to obtain $F_1,F_2$ clopens of $X$ such that
\begin{align*}
& F_1 \cup F_2 =A_2, \quad F_1 \cap F_2=\emptyset, \quad R[F_1] \subseteq A_2, \quad R[F_1] \nsubseteq F_1.
\end{align*}
It follows that $F_1 \cap R[F_2] \neq \emptyset$. Define $A_1'=A_1 \cup F_1$ and $A_2'=F_2$. Consequently, $A_1',A_2'$ are nonempty clopens such that for each $C \in \Pcal \setminus \{A\}$
\[
A_1' \cap R[C] \neq \emptyset \ \mbox{iff} \ C \in \S_1, \qquad A_2' \cap R[C] \neq \emptyset \ \mbox{iff} \ C \in \S_2 \qquad \mbox{and} \qquad A_1' \cap R[A_2'] \neq \emptyset.
\]
\end{proof}

Since the conditions (S1), (S2), (S3) correspond dually to the axioms (s1), (s2), (s3). Theorem~\ref{thm:axioms_model_completion} is an immediate consequence of Lemmas~\ref{lem:exclosedimpliesS1S2S3} and~\ref{lem:S1S2S3implyexclosed}.

We can use Theorem~\ref{thm:axioms_model_completion} to give another proof of the fact that the formula~(\ref{eq:ex}) corresponding to the $\Pi_2$-rule $(\rho_9)$ is provable in $\Con^\star$, hence 
$(\rho_9)$ is admissible in $\sf S^2IC$. 

\begin{corollary}
The formula
\begin{equation*}
x \prec y \Rightarrow \exists z ~(z \prec z \fowedge x \prec z \prec y)
\end{equation*}
is provable in $\Con^\star$.
\end{corollary}

\begin{proof}
Let $A$ be an existentially closed contact algebra and $a,b \in A$ such that $a \prec b$. We can assume that $a \neq b$, otherwise the claim is trivial. We have 
\begin{equation*}
(a \vee \neg b ) \wedge (b \wedge \neg a)=0 \quad  \mbox{and} \quad b \wedge \neg a \prec 1 = (b \wedge \neg a) \vee a \vee \neg b.
\end{equation*}
Thus, by the axiom (s1) applied to $b \wedge \neg a$, $a$, $\neg b$, there are $c_1,c_2 \in A$ such that 
\begin{align*}
c_1 \vee c_2=b \wedge \neg a, \quad c_1 \wedge c_2=0, \quad c_1 \neq 0, \quad c_2 \neq 0, \quad c_1 \prec c_1 \vee a, \quad c_2 \prec c_2 \vee \neg b.
\end{align*}
Since $a \prec b$ and $c_1 \prec c_1 \vee a \le b$, it follows that $a \vee c_1 \prec b$.
Moreover,
\begin{equation*}
a \vee c_1 = a \vee ((b \wedge \neg a) \wedge \neg c_2) = b \wedge \neg c_2= \neg(c_2 \vee \neg b) \prec \neg c_2
\end{equation*}
where the second equality is a consequence of $a \le b$ and $c_2 \le \neg a$.
Therefore, 
$a \vee c_1 \prec b \wedge \neg c_2= a \vee c_1$. 
Let $d=a \vee c_1$. Then $d \prec d$ and  $a \le d \le b$, which imply $a \prec d \prec b$.
\end{proof}

\begin{definition}\label{def: 7.2}
Let $\cS$ be a modal system. We say that a set of $\Pi_2$-rules $\Theta$ \emph{derives} a $\Pi_2$-rule $\rho$ if 
\[
T_{\cS} \cup \{\Pi(\theta) \mid \theta \in \Theta\} \vDash \Pi(\rho).
\]
We say that a set of admissible rules $\Theta$ is a \emph{basis of admissible rules for $\cS$} if it is a minimal set of rules that derives every admissible rule.
\end{definition}

\begin{theorem}
A basis of admissible rules for $\sf S^2IC$ is given by the following three rules.
\begin{align*}
(\rho_{s1}) \quad &\inference{[\forall]((p_1 \lor p_2 \leftrightarrow \varphi_1) \land \neg(p_1 \land p_2) \land \langle \exists \rangle p_1 \land \langle \exists \rangle p_2 \land (p_1 \simpa p_1 \lor \varphi_2) \land (p_2 \simpa p_2 \lor \varphi_3)) \to \chi}{[\forall](\langle \exists \rangle \varphi_1 \land \neg (\varphi_1 \land (\varphi_2 \lor \varphi_3)) \land (\varphi_1 \simpa \varphi_1 \lor \varphi_2 \lor \varphi_3)) \to \chi}\\[0.5cm]
(\rho_{s2}) \quad &\inference{[\forall]((p_1 \lor p_2 \leftrightarrow \varphi_1) \land \neg(p_1 \land p_2) \land \neg(p_1 \simpa \neg \varphi_2)\land \neg(p_2 \simpa \neg \varphi_2) \land (p_1 \simpa \neg p_2)) \to \chi}{[\forall](\neg(\varphi_1 \land \varphi_2) \land \neg (\varphi \simpa \neg \varphi_2)) \to \chi}\\[0.5cm]
(\rho_{s3}) \quad &\inference{[\forall](((p_1 \lor p_2)\to \varphi) \land \neg (p_1 \land p_2) \land (p_1 \simpa \varphi) \land \neg (p_1 \simpa p_2)) \to \chi}{\langle \exists \rangle \varphi \to \chi}
\end{align*}
where $\langle \exists \rangle:=\neg [\forall] \neg$.
\end{theorem}

\begin{proof}
If $\cS$ is a modal system with universal modality $[\forall]$, then in $T_{\cS}$ the formula $x \neq 0$ is equivalent to $\langle \exists \rangle x=1$. Moreover, if $t_1,t_2$ are terms, then in $T_{\cS}$ the first order formula
\[
\forall \ux (t_1(\ux)=1 \Rightarrow \exists \underline{y} (t_2(\ux,\underline{y})=1))
\]
is equivalent to
\[
\forall \ux, z ([\forall]t_1(\ux) \nleq z \Rightarrow \exists \underline{y} ([\forall]t_2(\ux,\underline{y}) \nleq z)).
\]
It is then straightforward to see that the axioms $(s1)$, $(s2)$, and $(s3)$ of Theorem~\ref{thm:axioms_model_completion} are equivalent to $\Pi(\rho_{s1})$, $\Pi(\rho_{s2})$, and $\Pi(\rho_{s3})$  in the theory of contact algebras (thought of as simple $\sf S^2I$-algebras). 
Thus, by Theorems~\ref{thm:mainadm} and~\ref{thm:axioms_model_completion}, if $\rho$ is an admissible rule, then it would be a consequence of $\Pi(\rho_{s1})$, $\Pi(\rho_{s2})$, and $\Pi(\rho_{s3})$. This implies that $(\rho_{s1})$, $(\rho_{s2})$, and $(\rho_{s3})$ form a basis of admissible rules for $\sf S^2IC$.
\end{proof}

\section{Conclusions}

In this paper we studied admissibility of $\Pi_2$-rules. We derived three strategies for recognizing admissibility for such rules. These strategies used interpolation, uniform interpolation, and model completions, respectively. We tested these methods on the symmetric strict implication calculus $\mathsf{S^2IC}$ and showed that admissibility of $\Pi_2$-rules is decidable in $\mathsf{S^2IC}$. We also proved that 
the model completion of the theory of contact algebras (simple algebraic models of $\mathsf{S^2IC}$) is finitely axiomatizable. This allowed us to show that there is a finite basis for admissible $\Pi_2$-rules in $\mathsf{S^2IC}$.
Below we discuss some potential directions for future work.
 
In the last part of the paper we showed that there is a finite basis of admissible $\Pi_2$-rules for $\mathsf{S^2IC}$. For this in Definition~\ref{def: 7.2},  for  a set of $\Pi_2$-rules  $\Theta$ and rule $\rho$ we defined when
$\Theta$ derives $\rho$. Namely, $\Theta$ derives $\rho$ if the first-order translation of rules in $\Theta$ derive the first-order translation of $\rho$. The definition of a basis of admissible $\Pi_2$-rules is based on this definition. 
We leave it as an open problem to define when a set of  $\Pi_2$-rule derives a $\Pi_2$-rule purely in terms of these rules without appealing to their first-order correspondents.

Another direction is to study connections with the literature on admissibility of standard inference rules in contact algebras~\cite{BG13}. Our non-standard  $\Pi_2$-rules have the particular shape outlined in 
Definition~\ref{defirules} and these trivialize if they are standard (i.e., if $p$ does not occur in the formula $F$ from the premise). However, it could be interesting to 
analyze more general formats for non-standard rules that  also  include standard inference rules.

Finally, it will be interesting to study extensions of our results to systems over a distributive lattice reduct. Among others it might be useful  to develop a framework encompassing the important (non-standard) density rule of fuzzy and 
many-valued 
logics \cite{MM07, TT84}.



\bibliographystyle{plain}
\bibliography{aiml20,mcmt}


\vspace{5em}

\noindent Nick Bezhanishvili: Institute for Logic, Language and Computation, University of Amsterdam, P.O. Box 94242, 1090 GE Amsterdam, The Netherlands.\\
\texttt{N.Bezhanishvili@uva.nl}

\bigskip 

\noindent Luca Carai: Universit\`a degli Studi di Salerno, Dipartimento di Matematica, Via Giovanni Paolo II, 132 - 84084 Fisciano (SA), Italy.\\
\texttt{lcarai@unisa.it} 

\bigskip 

\noindent Silvio Ghilardi:  Department of Mathematics, Universit\`a degli Studi di Milano, 20133 Milano, Italy, \\
 \texttt{silvio.ghilardi@unimi.it} 

\bigskip 

\noindent Lucia Landi: Department of Mathematics, Universit\`a degli Studi di Milano, 20133 Milano, Italy, \\
 \texttt{lucia.landi@studenti.unimi.it} 

\end{document}